\documentclass[11pt]{article}
\usepackage{amsmath,amssymb,amsthm,mathrsfs}
\usepackage[a4paper,margin=1in]{geometry}
\RequirePackage[numbers]{natbib}
\RequirePackage{url}
\RequirePackage{mathtools}
\mathtoolsset{showonlyrefs}

\newcommand{\E}{\mathbb{E}}
\newcommand{\Pbb}{\mathbb{P}}
\newcommand{\1}{\mathbf{1}}

\newcommand{\dd}{\mathrm{d}}
\newtheorem{theorem}{Theorem}
\newtheorem{proposition}{Proposition}
\newtheorem{lemma}{Lemma}
\newtheorem{corollary}{Corollary}
\newtheorem{remark}{Remark}

\newcommand{\jaro}[1]{({\bfseries (JS): #1})}

\title{From invasion to replacement: excursions of the Moran model}
\author{Jaromir Sant, Department of Statistics, Universidad Carlos III de Madrid\\
Ivana Valentić, Department of Mathematics, University of Zagreb}
\date{}

\begin{document}
\maketitle
\begin{abstract}
    An excursion-theoretic description of the continuous-time Moran model with recurrent mutation is developed. Exploiting its birth--death structure, the killed Moran process is characterized through its Green function, resolvent, Laplace transforms of hitting times, and spectral decomposition. A duality between the Moran model and a pure-death process is also established, yielding a spectral representation of its transition probabilities. The characterization of the killed process leads to a boundary-switching excursion representation in which returning excursions correspond to unsuccessful invasions of an allelic type, whilst switching excursions describe population-wide allelic replacement. The resulting regenerative structure provides the exact joint distribution and Laplace transforms for several path statistics, including the number, height and duration of unsuccessful invasions preceding replacement. Using Doob $h$-transforms, convergence of the rescaled Moran excursion measures and associated point processes, restricted to excursions exceeding a fixed positive height, to their Wright–Fisher counterparts is established. These transformations also show how the law of a Wright--Fisher diffusion killed upon hitting the boundary can be related to that of a Wright--Fisher diffusion with appropriately transformed mutation parameters. 
\end{abstract}

	\section{Introduction}
	The Moran model is a classical finite-population model for the stochastic evolution of allele frequencies \cite{Etheridge2012,Ewens,Moran1958}. In the two-allele setting considered here, it is a continuous-time birth--death process on $\{0,\ldots,N\}$, recording the number of individuals carrying a given allele in a population of fixed size $N$. The boundary states 0 and $N$ correspond to the loss and population-wide replacement of that allele, respectively, while recurrent mutation allows the process to leave either boundary. Under appropriate rescaling, the allele-frequency process converges to the Wright--Fisher diffusion on ([0,1]) \cite[Chapter 4, Section 8]{EthierKurtz}. 
	
	Excursion theory provides a framework for decomposing Markov processes into excursions away from a distinguished state or set \cite{Ito,PitmanYor2003}. For the Wright--Fisher diffusion, this approach was developed in \cite{JenkinsEtAlWF}, where a boundary-switching excursion representation was constructed and used to investigate path properties, including Hausdorff dimension and the mutual singularity of path measures corresponding to different mutation parameters. In the finite-state Moran setting, the construction is more direct: excursions from either boundary begin with a jump to the adjacent interior state, and successive excursions are separated by exponential holding times at the boundary. The entrance law is therefore explicit, and the excursion measure is determined by the law of the Moran process started from the adjacent state and killed upon reaching either boundary.
	
	In this paper, we develop an explicit excursion theory for the continuous-time Moran model with recurrent mutation and exploit its regenerative structure to obtain joint path statistics beyond conventional first-passage quantities \cite{Huillet2013}. Excursions from a boundary can be interpreted as successive invasion attempts: those returning to the starting boundary represent unsuccessful invasions, while an excursion reaching the opposite boundary completes a population-wide replacement. This yields a representation of the Moran process as a sequence of unsuccessful invasions followed by a successful one, with the construction restarting from the opposite boundary after each replacement. Using explicit Green functions, resolvent kernels and discounted hitting-time transforms, we obtain exact distributions and joint transforms for the number, heights and durations of unsuccessful invasions preceding replacement, together with the duration of the eventual successful invasion. The excursion framework thus retains information about the sequence of invasion attempts leading to replacement that is not captured by first-passage times alone.
	
	A key ingredient in the excursion construction is the analysis of the Moran process killed upon reaching prescribed boundary states, for which we obtain explicit Green functions, resolvent kernels and spectral representations. We further establish convergence of the killed Moran process and its associated excursion quantities to their Wright--Fisher diffusion counterparts under appropriate scaling. This involves a Doob $h$-transform for the killed processes: while the transformed Wright--Fisher diffusion remains a Wright--Fisher diffusion with modified mutation parameters, the transformed Moran chain does not retain the original Moran transition structure, thus requiring a separate convergence argument. The same transformation also yields a representation of the killed Wright--Fisher transition density in terms of the transition density of a Wright--Fisher diffusion with transformed mutation parameters, allowing the use of existing spectral expansions and exact simulation methods to evaluate quantities associated with the killed process \cite{JenkinsSpano,Santetal2026}. Finally, we establish a duality between the Moran model and a finite pure-death process, analogous to the classical duality between the Wright--Fisher diffusion and the coalescent \cite{Griffiths79,Tavare84,JenkinsSpano}, yielding an alternative representation of its transition probabilities through hypergeometric and beta--binomial sampling kernels.
	
	The paper is organized as follows. Section~\ref{sec:MoranModel} recalls the basic properties of the continuous-time Moran model, including its spectral decomposition, duality with a pure-death process, hitting-time transforms, and Green and resolvent kernels. Section~\ref{sec:KilledMoranModel} develops the corresponding results for the Moran process killed upon reaching two prescribed states, including the spectral decomposition of the killed generator. Section~\ref{sec:Excursions} constructs the boundary-switching excursion process through its local times, entrance laws and killed semigroups, while Section~\ref{sec:JointPathStatistics} exploits this construction to derive joint path statistics of successive invasion attempts. Section~\ref{sec:ConvergenceToDiffusion} establishes convergence to the corresponding Wright--Fisher diffusion quantities and obtains a representation of the killed Wright--Fisher transition density. Proofs of auxiliary results are collected in the Appendix \ref{appendix}.
	
    \section{The continuous-time neutral Moran model with mutation}\label{sec:MoranModel}
	Consider a population of fixed size $N$ in which each individual carries one of the two allelic types $a$ or $A$. Without loss of generality, let $X^N=(X_t^N)_{t\geq0}$ be the process that counts the number of individuals carrying allele $A$, such that $X^N$ takes values in $I_N=\{0,1,\ldots,N\}$. 
    
    Reproductive events occur at a constant rate $1/2$, and at each reproductive event, one individual is chosen uniformly at random to reproduce and independently another is chosen uniformly at random to be replaced. Note in particular that the same individual may be chosen for both roles. The offspring inherits the type of its parent unless a mutation occurs. Conditional on a reproductive event, an $A$-type offspring mutates to type $a$ with probability $u$, while an $a$-type offspring mutates to type $A$ with probability $v$. We shall assume that $v = \theta_1 N^{-1}$ and $u = \theta_2 N^{-1}$ for some $\theta_1, \theta_2 > 0$, and we define $\theta := \theta_1 + \theta_2$. 
    
    $X^N$ is a continuous-time birth--death chain with generator $Q_N$ acting on functions $f:I_N\to\mathbb{R}$ such that
	\begin{align}\label{eq:operator L^N}
		(Q_Nf)(i)&=\lambda_i\big(f(i+1)-f(i)\big)+\mu_i\big(f(i-1)-f(i)\big),
		\qquad 0\le i\le N,
	\end{align}
	where 
	\begin{align}\label{eq:lambdaMu}
		\lambda_i
		&=\frac{1}{2}\delta_N\left(1-\frac{i}{N}\right)\left(\frac{i}{N}+\frac{a_N}{N}\right) & 
		\mu_i
		&=\frac{1}{2}\delta_N\frac{i}{N}\left(1-\frac{i}{N}+\frac{b_N}{N}\right),
	\end{align}
	where we have defined 
	\begin{align}\label{eq:a_N_b_N}
		\delta_N &:= 1-\theta/N, & a_N &= \frac{\theta_1}{\delta_N}, & b_N &= \frac{\theta_2}{\delta_N}.
	\end{align}
	Throughout, we shall assume that $N>\theta$, so that $\delta_N>0$ and hence $a_N,b_N>0$. This is the classical continuous-time Moran model with mutation; see \cite{Moran1958,Moran1962,KarlinMcGregorMoran,Zhou2008,Etheridge2012}.
	
	Note that since $\theta_1,\theta_2>0$, the chain $X^N$ is irreducible on $I_N$ and therefore admits a unique stationary distribution. The following classical result (see for instance \cite{Moran1958,Ewens}) identifies it explicitly; we include a short proof for completeness. 
	\begin{proposition}\label{prop:stat_dis}
		The stationary law $\boldsymbol\pi := (\pi_i)_{i\in I_N}$ of $X^N$ is the $\text{BetaBin}(N,a_N,b_N)$ distribution, where $a_N$, and $b_N$ are defined by \eqref{eq:a_N_b_N}.
	\end{proposition}
	\begin{proof}
		Using standard birth-and-death process theory, set
		\begin{equation*}
			\tilde{\pi}_0:=1,\qquad
			\tilde{\pi}_i:=\prod_{k=0}^{i-1}\frac{\lambda_k}{\mu_{k+1}},\qquad
			\tilde{\Pi}:=\sum_{j=0}^{N}\tilde{\pi}_j,\qquad
			\text{so }\ \ \pi_i=\tilde{\pi}_i/\tilde{\Pi}.
		\end{equation*}
		By rearranging we get 
		\begin{align*}
			\tilde{\pi}_i&= \prod_{k=0}^{i-1}\frac{\lambda_k}{\mu_{i-k}}=\prod_{k=0}^{i-1} \frac{N-k}{i-k}\prod_{k=0}^{i-1} \frac{a_N+k}{N-i+b_N+k}=\binom{N}{i} \frac{B\left(a_N+i, N+b_N-i\right)}{B\left(a_N,b_N+N\right)},
		\end{align*}
		from which it follows that
		\begin{align*}
			\tilde{\Pi}:=&\sum_{j=0}^{N} \binom{N}{j} \frac{B\left(a_N+j,N+b_N-j\right)}{B\left(a_N,b_N+N\right)}=\frac{B\left(a_N,b_N\right)}{B\left(a_N,N+b_N\right)}.
		\end{align*}
		Thus, for $I_N \sim \text{BetaBin}(N,a_N,b_N)$, 
		\begin{equation*}
			\Pbb(I_N=i)=\binom{N}{i}\frac{B(a_N+i,N+b_N-i)}{B(a_N,b_N)}=\frac{\tilde{\pi}_i}{\tilde{\Pi}}=\pi_i.
		\end{equation*}
	\end{proof}
    Observe that the scale function of the process $X^N$ is given by
    \begin{align}\label{eq:scaleFn}
        s_0^{(N)}=0, \,\,\, s^{(N)}_n := \sum_{k=0}^{n-1} \frac{1}{\pi_k \lambda_k} \qquad n \in \{1,2,\ldots,N\}.
    \end{align}
    
	Since $X^N$ is a birth--death process, the stationary distribution $\boldsymbol\pi$ satisfies the detailed-balance equation $\pi_i\lambda_i=\pi_{i+1}\mu_{i+1}$, for $0\leq i<N$, and hence $X^N$ is reversible with respect to $\boldsymbol\pi$. Consequently, $Q_N$ is self-adjoint on $\ell^2(\boldsymbol\pi)$ and admits an orthogonal basis of eigenfunctions \citep{KarlinMcGregor}. In what follows, let $\mathbb{P}_{i}^{N}$ denote the law of the process when started from $i\in I_N$ and denote by $\mathbb{E}_{i}^{N}$ the expectation operator with repect to this probability.
    
    In light of the above, the continuous time Moran model with mutation has transition probabilities given by \citep{KarlinMcGregorMoran,Dominguez2021}
	\begin{equation}\label{eq:transitionProbabilityHahn}
		\mathbb{P}_i^N(X_t^N=j) = \pi_j \sum_{n=0}^{N} e^{-\beta_n t} \frac{H_n(i;a_N,b_N,N)H_n(j;a_N,b_N,N)}{\langle H_n(\cdot;a_N,b_N,N),H_n(\cdot;a_N,b_N,N) \rangle_{\boldsymbol\pi}},
	\end{equation}
	where for $\rho \in \mathbb{C}$ we define
    \begin{equation}\label{beta_rho}
        \beta_\rho := \frac{\delta_N}{2N^2} \rho(\rho+a_N+b_N-1),
    \end{equation}
    for $f, g : I_N \mapsto \mathbb{R}$ we set $\langle f,g \rangle_{\boldsymbol\pi} := \sum_{i=0}^N f(i)g(i)\pi_i$, and $H_n(i;a_N,b_N,N)$ denotes the Hahn polynomials defined as
	\begin{align}\label{eq:HahnPolynomial}
		H_n(x;a_N,b_N,N) &={}_3F_2 \left( \begin{matrix}
												-n,\ n+a_N+b_N-1,\ -x\\
												a_N,\ -N
										\end{matrix};1\right) \\
		&=\sum_{k=0}^n \frac{(-n)_k(n+a_N+b_N-1)_k(-x)_k}{(a_N)_k(-N)_k k!},
	\end{align} 
	for $n,x \in \{0, \dots, N\}$, and where $(z)_k=z(z+1)\cdots(z+k-1)$ denotes the rising factorial, with $(z)_0 := 1$. For notational convenience, we will henceforth suppress the explicit dependence of $H_n(\cdot; a_N, b_N, N)$ on $a_N, b_N$ and $N$ and simply write $H_n(\cdot)$. The Hahn polynomials are orthogonal with respect to the above defined scalar product, have squared norm
    \begin{equation}\label{eq:HahnNorm}
         \langle H_n,H_n\rangle_{\boldsymbol\pi}=
         \begin{cases}
         1,&n=0,\\[2pt]
         \displaystyle\frac{(N+a_N+b_N)_n(b_N)_n}
         {\binom Nn(2n+a_N+b_N-1)(a_N+b_N)_{n-1}(a_N)_n},&1\leq n\leq N.
         \end{cases}
    \end{equation}
	and $H_n(0)=1$. For more details on Hahn polynomials, see \cite[Section~1.5]{KoekoekSwarttouw1996}.

    \subsection{A dual process and transition probability representation}\label{sec:Duality}
	A classical result in mathematical population genetics is the moment duality between the Wright--Fisher diffusion and the Kingman coalescent \citep{EthierGriffiths93,GriffithsSpano}, which in turn has allowed for an explicit decomposition of the diffusion's transition density in terms of the the transition probabilities of the coalescent \citep{Griffiths80,Tavare84}. In this section we show how the continuous-time Moran model is related to a pure-death process under moment duality, and exploit this to get an explicit expression for the transition probabilities of the former. 
	
	For $0 \leq i \leq N$ and $0\leq m\leq N$, let
	\begin{equation}\label{eq:MoranDualHypergeometricKernel}
		\mathcal H_{N,m;i}(\ell) :=\frac{\binom{i}{\ell}\binom{N-i}{m-\ell}}{\binom{N}{m}}, \qquad (0 \vee m+i-N) \leq \ell\leq (i \wedge m).
	\end{equation}
	denote the hypergeometric probability mass function. Given $\ell$, define for $\ell \leq j \leq N - m + \ell$ 
	\begin{equation}\label{eq:MoranDualPosteriorKernel}
		\mathcal D_{N;m,\ell}(j) :=\binom{N-m}{j-\ell} \frac{B(a_N+j,b_N+N-j)}{B(a_N+\ell,b_N+m-\ell)},
	\end{equation}
	and set $\mathcal D_{N;m,\ell}(j) = 0$ for any other $j$.  Set
	\begin{equation}\label{eq:MoranDualLinkKernel}
		K_{m}^{N}(i,j):=\sum_{\ell=0}^{m} \mathcal H_{N,m;i}(\ell)\mathcal D_{N;m,\ell}(j),
	\end{equation}
	and observe that $K_{m}^N(i,j)$ depends on $i$ only through $\mathcal{H}_{N,m;i}(\ell)$. Thus, for $i \in I_N$, the action of $Q_N$ on $K_m^N(\cdot,j)$ is given by
	\begin{equation}
		Q_N K_{m}^{N}(i,j) = \lambda_i\left(K_{m}^{N}(i+1,j) - K_{m}^{N}(i,j)\right) + \mu_i\left(K_{m}^{N}(i-1,j) - K_{m}^{N}(i,j)\right).
	\end{equation}
	The following lemma (whose proof we defer to Appendix \ref{app:Sec2Proofs}) will be useful in what follows 
    \begin{lemma}\label{lem:gammam}
		Define $\gamma_{m}^{N} := \frac{\delta_{N}}{2N^2}m(m+a_N + b_N -1)$. Then
		\begin{align}
			Q_N \mathcal{H}_{N,m;i}(\ell) &= -\gamma_{m}^{N}\mathcal{H}_{N,m;i}(\ell) + \frac{\delta_{N}}{2N^2} m(a_N+\ell-1)\mathcal{H}_{N,m-1;i}(\ell-1) \\
			&{}+ \frac{\delta_{N}}{2N^2}m(b_N + m - \ell - 1)\mathcal{H}_{N,m-1;i}(\ell).
		\end{align}
	\end{lemma}
    
    Lemma \ref{lem:gammam} gives immediately that
	\begin{align}
		Q_N K_{m}^{N}(i,j) &= -\gamma_{m}^{N}K_{m}^{N}(i,j) + \frac{\delta_N}{2N^2} m \sum_{\ell=0}^{m-1}\mathcal{H}_{N,m-1;i}(\ell)\Big[(a_N+\ell)\mathcal{D}_{N;m,\ell+1}(j) \nonumber \\
        &\qquad\qquad\qquad\qquad\qquad\qquad\qquad\qquad{}+ (m+b_N-\ell-1)\mathcal{D}_{N;m,\ell}(j)\Big] \nonumber \\
        &= \gamma_{m}^{N}\left(K_{m-1}^{N}(i,j) - K_{m}^{N}(i,j)\right),
	\end{align}
	where in the first equality we have shifted indices and collected like terms, and in the second equality we use Pascal's identity. Since $\gamma_{m}^{N} \ge 0$ for any $m \ge 0$, we recognize the RHS as the generator of a pure-death process, giving the following result.
	
	\begin{theorem}\label{thm:DualityMoranPureDeath}
		Let $(X^{N}_{t})_{t\geq0}$ be a continuous time Moran model with generator $Q_N$ given by \eqref{eq:operator L^N}. Let $(M_{t}^{N})_{t\geq0}$ be a pure-death process with transition rates $\gamma_m^N=\frac{\delta_N}{2N^2}m(m+a_N+b_N-1)$ for $m \in \{ 0, \dots, N \}$. Then for any $t \geq 0$
        \begin{align}\label{eq:MoranDuality}
			\mathbb{E}\left[K_{m}^{N}(X_{t}^{N}, j) | X_{0}^{N} = x\right] = \mathbb{E}\left[K_{M_{t}^{N}}^N(x,j) | M_{0}^{N} = m\right].
		\end{align}
	\end{theorem}
    
	Denote by $\mathcal{B}_N$ the generator of the above pure-death process $(M_t^{N})_{t\geq0}$. We point out that the domain of either generator is $\mathcal{D}(Q_N) = \mathcal{D}(\mathcal{B}_N) = \mathbb{R}^{I_N}$. Using \eqref{eq:MoranDuality}, we can derive an explicit representation for the transition probability of the process $(X_{t}^{N})_{t\geq0}$. 

    \begin{theorem}\label{thm:MoranTransitionProbability}
        Let $(X_t^N)_{t\ge0}$ and $(M_t^N)_{t\ge0}$ be as in Theorem \ref{thm:DualityMoranPureDeath}. Then for any $i, j \in I_N$ and $t \ge 0$ we have that
        \begin{align}\label{eq:transitionDensityDecompositionFull}
    		\mathbb{P}[X_{t}^{N} = j | X_{0}^{N} = i] &= \sum_{m=0}^{N}q_{N,m}^{N}(t)\sum_{\ell=0}^{m}\mathcal{H}_{N,m;i}(\ell)\mathcal{D}_{N;m,\ell}(j),
    	\end{align}
        where
        \begin{align}\label{eq:qm(t)Full}
    		q_{N,m}^{N}(t) &= \sum_{r=m}^{N}(-1)^{r-m}\frac{N!}{m!(r-m)!(N-r)!}\frac{2r + a_N + b_N - 1}{m + r + a_N + b_N - 1} \nonumber\\
    		&\qquad\qquad{}\times\frac{(m+a_N+b_N)_{(N-m)}}{(m+r+a_N+b_N)_{(N-m)}}e^{-\frac{\delta_{N}}{2N^2}r(r+a_N+b_N-1)t}.
    	\end{align}
    \end{theorem}
    \begin{proof}
        Let $P_{t}^{N}$ denote the semigroup operator associated with $(X_{t}^{N})_{t\geq0}$ and $Q_{t}^{N}$ that associated with the pure-death process $(M_t^{N})_{t\geq0}$. Write
	$q_{n,m}^{N}(t):=\mathbb{P}[M_t^{N}=m \mid M_{0}^{N} = n]$ and $p_{i,j}^{N}(t) := \mathbb{P}[X_t^{N} = j \mid X_0^{N} = i]$. Then
	\begin{align}\label{eq:transitionDensityDecomposition}
		p_{i,j}^{N}(t) = \mathbb{P}[X_{t}^{N} = j | X_{0}^{N} = i] = \left(P_{t}^{N} K_{N}^{N}(\cdot,j)\right)(i) &= \sum_{k=0}^{N}p_{i,k}^{N}(t)K_{N}^{N}(k,j) \nonumber \\
		&= \sum_{m=0}^{N}q_{N,m}^{N}(t)K_{m}^{N}(i,j),
	\end{align}
	where we have used Theorem \ref{thm:DualityMoranPureDeath} and $K_{N}^{N}(i,j) = \delta_{ij}$. Then \eqref{eq:transitionDensityDecompositionFull} follows immediately from \eqref{eq:MoranDualLinkKernel}. It remains to prove \eqref{eq:qm(t)Full}. In particular the forward equation for the pure-death process is given by 
	\begin{align}\label{eq:forwardEqnDual}
		\frac{d}{dt}q_{N,m}^{N}(t) = \gamma_{m+1}^{N}q_{N,m+1}^{N}(t) - \gamma_{m}^{N}q_{N,m}^{N}(t) & \textnormal{ for } m < N
	\end{align}
	with boundary conditions $q_{N,N}(0) = 1$ and $q_{N,m}(0) = 0$ for $m < N$. By induction, $q_{N,m}^{N}(t)$ is given by 
	\begin{align}\label{eq:InductiveHypothesis}
		q_{N,m}^{N}(t) &= \left(\prod_{\ell = m+1}^{N}\gamma_{\ell}^{N}\right)\sum_{r=m}^{N}\frac{e^{-\gamma_r^{N}t}}{\prod_{\substack{{\ell = m} \\ {\ell \neq r}}}^{N}(\gamma_{\ell}^{N} - \gamma_{r}^{N})}.
	\end{align}
	Indeed when $m=N$, it is immediate that $q_{N,N}(t) = e^{-\gamma_{N}t}$. Now suppose \eqref{eq:InductiveHypothesis} holds for $q_{N,m+1}^{N}(t)$, we solve \eqref{eq:forwardEqnDual} to get
	\begin{align*}
		q_{N,m}^{N}(t) &= \left(\prod_{\ell = m+1}^{N}\gamma_{\ell}^{N}\right) \sum_{r=m+1}^{N}\frac{e^{-\gamma_{m}^{N}t} - e^{-\gamma_{r}^{N}t}}{(\gamma_{r}^{N} - \gamma_{m}^{N})\prod_{\substack{{\ell = m+1} \\ {\ell \neq r}}}^{N}(\gamma_{\ell}^{N} - \gamma_{r}^{N})}.
	\end{align*}
    Standard partial fraction identities allow us to deduce that
	\begin{align*}
		\sum_{r=m+1}^N\frac{1}{(\gamma_r^N-\gamma_m^N)
       \prod_{\substack{\ell=m+1\\\ell\ne r}}^N(\gamma_\ell^N-\gamma_r^N)} =\frac1{\prod_{\ell=m+1}^N(\gamma_\ell^N-\gamma_m^N)},
	\end{align*}
	from which \eqref{eq:InductiveHypothesis} follows. The products of $\gamma_{m}^{N}$ terms can be further expanded to get
	\begin{align*}
		\prod_{\ell = m+1}^{N}\gamma_{\ell}^{N} = \left(\frac{\delta_{N}}{2N^2}\right)^{N-m}\frac{N!}{m!}\frac{\Gamma(N+a_N+b_N)}{\Gamma(m+a_N+b_N)} = \left(\frac{\delta_{N}}{2N^2}\right)^{N-m}\frac{N!}{m!}(m+a_{N}+b_{N})_{(N-m)},
	\end{align*}
	  and since $\gamma_{\ell}^{N} - \gamma_{r}^{N} = \frac{\delta_{N}}{2N^2}(\ell - r)(\ell + r + a_N + b_N -1)$ we further deduce that
	\begin{align*}
		\prod_{\substack{{\ell = m} \\ {\ell \neq r}}}^{N} \left(\gamma_{\ell}^N - \gamma_{r}^N\right) &= \left(\frac{\delta_N}{2N^2}\right)^{N-m}(-1)^{r-m}(r-m)!(N-r)!\frac{(m+r+a_N+b_N-1)_{(N-m+1)}}{2r + a_N + b_N -1}.
	\end{align*}
	from which \eqref{eq:qm(t)Full} follows.    
    \end{proof}
	
	The above representations \eqref{eq:transitionDensityDecompositionFull} and \eqref{eq:qm(t)Full} provide a finite-population counterpart of the classical Wright--Fisher transition density decomposition; see \cite{Griffiths79,Tavare84}, as well as \cite[equations (4) and (5)]{JenkinsSpano}. Thus Theorems \ref{thm:DualityMoranPureDeath} and \ref{thm:MoranTransitionProbability} give an exact discrete analogue of the ancestral mixture structure underlying the diffusion transition density. The explicit transition probabilities \eqref{eq:qm(t)Full} of the finite ancestral death process are themselves known; see \cite[equation (2.4)]{GriffithsSpano}. 
	
    \subsection{The Laplace transform of the hitting time}\label{sec:MoranModel_LT}
    Let $k\in I_N$, and define the hitting time $\tau_k^N := \inf\{t\geq0:X_t^N=k\}$. The standard formula for a birth--death chain gives
		\begin{align}\label{eq:hittingProbabilitySwitching}
			\mathbb P_i^N(\tau_k^N<\tau_j^N)
			&=\frac{s^{(N)}_i-s^{(N)}_j}{s^{(N)}_k-s^{(N)}_j}.
		\end{align}
    
    Observe that the Laplace transform of $\tau_k^N$ for any $\psi \geq 0$ is given by \cite[equation (9)]{PitmanYor2003} 
    \begin{align}\label{eq:LTHittingTime}
        \mathbb{E}_{x}^{N}\left[ e^{-\psi \tau_{k}^{N}} \right] = \begin{cases}
            \frac{u_{\psi}^{-}(x)}{u_{\psi}^{-}(k)} & x < k \\
            \frac{u^{-}_{\psi}(x)}{u^{+}_{\psi}(k)} & x > k,
        \end{cases}
    \end{align}
    where $u_{\psi}^{\pm}$ are non-negative solutions of $Q_N u_{\psi}^{\pm} = \psi u_{\psi}^{\pm}$ with appropriate boundary conditions, such that $u_{\psi}^{-}$ is increasing and $u_{\psi}^{+}$ is decreasing.

	Now let $\rho \in \mathbb{C}$ and $x\in I_N$, and define
	\begin{align}
		S_{\rho}^{-}(x) &:= \sum_{j=0}^{x} \frac{(-\rho)_{j}(\rho + a_N + b_N - 1)_{j}(-x)_{j}}{(a_N)_{j}(-N)_{j}j!}, \label{eq:S-rho-minus}\\
		S_{\rho}^{+}(x) &:= \sum_{j=0}^{N-y}\frac{(-\rho)_{j}(\rho + a_N + b_N - 1)_{j}(-(N-x))_{j}}{(b_N)_{j}(-N)_{j}j!}, \label{eq:S-rho-plus}
	\end{align}
	Then $S_{\rho}^{-}(x)$ and $S_{\rho}^{+}(x)$ respectively solve
	\begin{align}
		Q_N S_{\rho}^{-}(x) &= -\beta_{\rho}S_{\rho}^{-}(x), & \text{ for } 0 \leq &x \leq N-1, \\
		Q_N S_{\rho}^{+}(x) &= -\beta_{\rho}S_{\rho}^{+}(x), & \text{ for } 1 \leq &x \leq N.
	\end{align}
	where $\beta_\rho$ is given by \eqref{beta_rho}. 
	
	  For $x < k$ apply the generator $Q_N$ to $S_{\rho}^{-}(x) / S_{\rho}^{-}(k)$ for some $\rho \in \mathbb{C}$ to be determined, whilst for $x > k$ apply $Q_N$ to $S_{\rho}^{+}(x) / S_{\rho}^{+}(k)$. Either case leads to
	\begin{align}
		\lambda_x \frac{S_{\rho}^{\pm}(x+1) - S_{\rho}^{\pm}(x)}{S_{\rho}^{\pm}(k)} + \mu_x \frac{S_{\rho}^{\pm}(x-1) - S_{\rho}^{\pm}(x)}{S_{\rho}^{\pm}(k)} &= -\frac{\delta_N}{2N^2} \rho(\rho+a_N+b_N-1)u_x.
	\end{align}
	Equating $\psi=-\frac{\delta_N}{2N^2}\rho(\rho+a_N+b_N-1)$ and solving for $\rho$ we get
	\begin{equation}\label{eq:rhoSoln}
		\rho_N(\psi) := -\frac{a_N+b_N-1}{2} \pm \frac{\sqrt{(a_N+b_N-1)^2-\frac{8N^2\psi}{\delta_N}}}{2},
	\end{equation}
	where we have now made explicit the dependence of $\rho$ on both $N$ and $\psi$. Note that either choice in \eqref{eq:rhoSoln} leads to the same functions $S_\rho^\pm(\cdot)$, since the latter is symmetric for this choice of parameters. Thus we deduce that the functions in the RHS of \eqref{eq:LTHittingTime} are given by
    \begin{align}\label{eq:LTHittingTimesLeftRight}
        u_{\psi}^{-}(x) &= S_{\rho_N(\psi)}^{-}(x), & u_\psi^{+}(x) &= S_{\rho_N(\psi)}^{+}(x).
    \end{align}
	
	\subsection{The Casoratian}\label{sec:MoranModel_Casoratian}
	For functions $f,g:I_N\to\mathbb C$, set
	\begin{equation}\label{eq:discreteCasoratianDefinition}
		\mathcal W_i^N(f,g):=f(i)g(i+1)-f(i+1)g(i),
		\qquad 0\leq i\leq N-1,
	\end{equation}
	then for $\psi > 0$ the scale function normalized Casoratian for the continuous-time Moran model with mutation as
	\begin{equation}\label{eq:unrestrictedScaleCasoratian}
		w_{\psi}^N
		:=\frac{-\mathcal{W}^{N}_i
			(S_{\rho_N(\psi)}^{-},S_{\rho_N(\psi)}^{+})}{s^{(N)}_{i+1}-s^{(N)}_i},
	\end{equation}
	where we observe that the right-hand side does not depend on $i$. Abel's identity allows us to conclude that
	\begin{align}
		\mathcal{W}^{N}_i(S_{\rho_N(\psi)}^{-},S_{\rho_N(\psi)}^{+}) &= \frac{i(N-i+b_N)}{(N-i)(i+a_N)}\mathcal{W}^{N}_{i-1}(S_{\rho_N(\psi)}^{-},S_{\rho_N(\psi)}^{+}) \nonumber\\
		&= \mathcal{W}_0^{N} \Gamma(i+1)\frac{\Gamma(N+b_N)\Gamma(N-i)\Gamma(a_N+1)}{\Gamma(N)\Gamma(N-i+b_N)\Gamma(a_N+i+1)},
	\end{align}
	where for brevity we write 
	\begin{align*}
		\mathcal{W}^{N}_0 := \mathcal{W}^{N}_0(S_{\rho_N(\psi)}^{-},S_{\rho_N(\psi)}^{+})=S_{\rho_N(\psi)}^{-}(0)S_{\rho_N(\psi)}^{+}(1) - S_{\rho_N(\psi)}^{-}(1)S_{\rho_N(\psi)}^{+}(0).
	\end{align*} 
	
	\subsection{Green's function}\label{sec:MoranModel_Green}
	Denote by $_{\boldsymbol\pi}\mathbb{P}_{i}^{N}\left[X_t^N = j\right]$ the transition probability of the continuous-time neutral Moran model with respect to its stationary distribution $\boldsymbol\pi$. Then from \eqref{eq:transitionProbabilityHahn} we get that
	\begin{align}
		_{\boldsymbol\pi}\mathbb{P}_{i}^{N}\left[X_t^N = j\right] = \sum_{n=0}^{N}e^{-\beta_n t}\frac{H_{n}(i)H_{n}(j)}{\langle H_n, H_n\rangle_{\pi}}.
	\end{align}
	Next we define for $\psi > 0$ the Green's function as 
	\begin{align}\label{eq:GreenFnDefn}
		G_{\psi}^{N}(i,j) = \int_{0}^{\infty}e^{-\psi t}{}_{\boldsymbol\pi}\mathbb{P}_{i}^{N}\left[X_t^N = j\right]dt.
	\end{align}
	We have the following result for the Green's function
	\begin{proposition}\label{prop:GreenFns}
		The Green's function for the continuous time Moran model is given by
		\begin{align}
			G^{N}_{\psi}(i,j) &= \sum_{n=0}^{N}\frac{1}{\psi + \beta_n}\frac{H_{n}(i)H_{n}(j)}{\langle H_n, H_n\rangle_{\pi}} \label{eq:GreenFunctionsMoran1}\\ 
			&= \frac{S_{\rho_N(\psi)}^{-}(i\wedge j)S_{\rho_N(\psi)}^{+}(i \vee j)}{w_{\psi}^{N}} \label{eq:GreenFunctionsMoran2}\\ 
			&= \frac{1}{\pi_j}\sum_{m=0}^{N}\frac{1}{\psi+\beta_m}\prod_{r=m+1}^{N}\frac{\beta_r}{\beta_r + \psi}\sum_{\ell=0}^{m}\mathcal{H}_{N,m;i}(\ell)\mathcal{D}_{N;m,\ell}(j)\label{eq:GreenFunctionsMoran3}
		\end{align}
	\end{proposition}
	\begin{proof}
		Equation \eqref{eq:GreenFunctionsMoran1} follows by observing that
		\begin{align}
			G_{\psi}^{N}(i,j) = \int_{0}^{\infty}e^{-\psi t}\sum_{n=0}^{N}e^{-\beta_n t}\frac{H_{n}(i)H_{n}(j)}{\langle H_n, H_n\rangle_{\pi}}dt = \sum_{n=0}^{N}\frac{1}{\psi + \beta_n}\frac{H_{n}(i)H_{n}(j)}{\langle H_n, H_n\rangle_{\pi}},
		\end{align}
		whereas \eqref{eq:GreenFunctionsMoran2} follows immediately from \cite[equation (37)]{PitmanYor2003}. For \eqref{eq:GreenFunctionsMoran3} we make use of Theorem \ref{thm:MoranTransitionProbability} as well as the arguments in Proposition 2.2 in \cite{JenkinsEtAlWF} with the difference that here our products are finite.
	\end{proof}
	
	The resolvent is defined as 
	\begin{align}\label{eq:resolventDefn}
		R_{\psi}^{N}f(i) := \mathbb{E}_{i}^{N}\left[\int_{0}^{\infty}e^{-\psi t}f(X_t^N)dt\right] = \sum_{j=0}^{N}G_{\psi}^{N}(i,j)f(j)\pi_j,
	\end{align}
	which coupled with Proposition \ref{prop:GreenFns} gives the following.
	\begin{proposition}
		For the continuous time Moran model, the resolvent is given by
		\begin{align}
			R_{\psi}^{N}f(i) &= \sum_{n=0}^{N}\frac{1}{\psi + \beta_n}\frac{H_{n}(i)}{\langle H_n, H_n\rangle_{\boldsymbol\pi}}\mathbb{E}_{\boldsymbol\pi}\left[H_n(J)f(J)\right] \label{eq:Resolvent1} \\
			&= \frac{1}{w_{\psi}^N}\mathbb{E}_{\boldsymbol\pi}\left[S_{\rho_N(\psi)}^{-}(i \wedge J)S_{\rho_N(\psi)}^{+}(i \vee J)f(J)\right] \label{eq:Resolvent2} \\
			&= \sum_{m=0}^{N}\frac{1}{\psi+\beta_m}\prod_{r=m+1}^{N}\frac{\beta_r}{\beta_r + \psi}\mathbb{E}\left[\sum_{j=0}^{N}\mathcal{D}_{N;m,K}(j)f(j)\right] \label{eq:Resolvent3}
		\end{align}
		where the expectation in \eqref{eq:Resolvent3} is taken with respect to $K \sim \textnormal{Hyp}(N,m,i)$.
	\end{proposition}
	\begin{proof}
		Follows directly from \eqref{eq:resolventDefn} and Proposition \ref{prop:GreenFns}.
	\end{proof}
	
	\section{Killed Moran model}\label{sec:KilledMoranModel}
    We now proceed to derive a number of properties of the continuous time Moran model killed when reaching points $a < b$ for $a, b \in I_N$. In particular, define $I_{N}^{[a,b]} := \{a,\ldots,b\}$, $I_{N}^{(a,b)} := \{ a+1, \dots, b-1 \}$, $\partial_{a,b}:=\{a,b\}$, set $\tau^{N}_{a,b}:=\tau_a^N\wedge\tau_b^N$, and let $Q_{N}^{(a,b)}$ denote the restriction of the generator $Q_N$ to $I_{N}^{(a,b)}$. Then for $i,j\in I_{N}^{(a,b)}$, we denote the transition probabilities of the continuous-time Moran model with generator \eqref{eq:operator L^N} killed when it first reaches $\partial_{a,b}$ by
	\begin{equation}
		p_{i,j}^{N,(a,b)}(t)
		:=\mathbb P_i^N(X_t^N=j,\ t<\tau^{N}_{a,b}).
		\label{eq:killedMoranTransitionProbability}
	\end{equation}
	We further let $P_t^{N, (a,b)}$ denote the corresponding semigroup operator, such that $P_t^{N,(a,b)}=\exp(tQ_N^{(a,b)})$. The rates $\mu_{a+1}$ and $\lambda_{b-1}$ retained on the diagonal represent killing at $a$ and $b$ respectively. For $\psi\geq 0$, similarly to \eqref{eq:GreenFnDefn} and \eqref{eq:resolventDefn}, define the killed Green's function as
    \begin{align}
		G_\psi^{N,(a,b)}(i,j) &:=\int_0^\infty e^{-\psi t}\frac{p_{i,j}^{N,(a,b)}(t)}{\pi_j}\,\dd t,
		\label{eq:killedMoranResolventGreen}
	\end{align}
	as well as the killed resolvent by
	\begin{align}
		R_\psi^{N,(a,b)}f(i) &:= \mathbb{E}_i^{N}\left[\int_{0}^{\infty}e^{-\psi t}f(X_t^{N})\mathbf1_{\{ t < \tau_{a,b}^{N} \}}dt\right] \nonumber \\
        &=\mathbb E_i^N\left[\int_0^{\tau_{a,b}^N}e^{-\psi t}f(X_t^N)\,\dd t \right] = \sum_{j=1}^{N-1}G_{\psi}^{N,(a,b)}(i,j)f(j)\pi_{j}.
		\label{eq:killedMoranResolventDefinition}
	\end{align}
	
	\subsection{Fundamental solutions and the killed resolvent}\label{sec:KilledMoranModel_FundamentalSoln}
	We now derive the explicit form of the killed Green's function and the killed resolvent for the continuous time Moran model. For $\psi > 0$, set
	\begin{align*}
		U_{a,b}^{+}(i; \psi) &:= \mathbb{E}_{i}\left[ e^{-\psi \tau_a^N} \mathbf1_{\{ \tau_{a}^{N} < \tau_{b}^{N} \}} \right] & U_{a,b}^{-}(i; \psi) &:= \mathbb{E}_{i}\left[ e^{-\psi \tau_b^N} \mathbf1_{\{ \tau_{b}^{N} < \tau_{a}^{N} \}} \right].
	\end{align*}
	We have the following lemma:
	\begin{lemma}\label{lem:Uab}
		For $\psi>0$, $\rho_N(\psi)$ as in \eqref{eq:rhoSoln}, and $i\in \{a, a+1, \dots, b-1, b\}$, 
		\begin{align}
			U^{+}_{a,b}(i; \psi) &= \frac{S^{-}_{\rho_N}(b)S^{+}_{\rho_N}(i) - S^{+}_{\rho_N}(b)S^{-}_{\rho_N}(i)}{S^{-}_{\rho_N}(b)S^{+}_{\rho_N}(a) - S^{+}_{\rho_N}(b)S^{-}_{\rho_N}(a)}, \label{eq:UabMinus}\\ 
			U^{-}_{a,b}(i; \psi) &= \frac{S^{-}_{\rho_N}(i)S^{+}_{\rho_N}(a) - S^{-}_{\rho_N}(a)S^{+}_{\rho_N}(i)}{S^{-}_{\rho_N}(b)S^{+}_{\rho_N}(a) - S^{+}_{\rho_N}(b)S^{-}_{\rho_N}(a)}. \label{eq:UabPlus}
		\end{align}
		
    In particular when $a = 0$ and $b = N$ we get that
    	\begin{align}
    		\Phi^{N}_{\psi,+}(i) := U^{+}_{0,N}(i; \psi) &= \frac{S_{\rho_N(\psi)}^{-}(i) - S_{\rho_N(\psi)}^{-}(N)S_{\rho_N(\psi)}^{+}(i)}{1-S_{\rho_N(\psi)}^{-}(N)S_{\rho_N(\psi)}^{+}(0)} \\
    		\Phi^{N}_{\psi,-}(i) := U^{-}_{0,N}(i; \psi) &= \frac{S_{\rho_N(\psi)}^{+}(i) - S_{\rho_N(\psi)}^{+}(0)S_{\rho_N(\psi)}^{-}(i)}{1-S_{\rho_N(\psi)}^{-}(N)S_{\rho_N(\psi)}^{+}(0)} 
    	\end{align}
    \end{lemma}
    \begin{proof}
    Follows by observing that $S_{\rho_N(\psi)}^{\pm}(\cdot)$ solve $Q_N S_{\rho_N(\psi)}^{\pm}(i) = -\beta_{\rho_N(\psi)}S_{\rho_N(\psi)}^{\pm}(i)$ for $i \in I_{N}{^{(a,b)}}$, and that $U^{+}_{a,b}(a; \psi) = 1$, $U^{+}_{a,b}(b; \psi) = 0$, $U^{-}_{a,b}(a; \psi) = 0$, and $U^{-}_{a,b}(b; \psi) = 1$.
	\end{proof}
	
	Fix $k\in I_{N}^{(a,b)}$ and define the corresponding Casoratian as
	\begin{equation}\label{eq:killedScaleCasoratian}
		w_\psi^{N,(a,b)}
		:=-\frac{
			U_{a,b}^{-}(k;\psi)U_{a,b}^{+}(k+1;\psi) - U_{a,b}^{-}(k+1;\psi)U_{a,b}^{+}(k;\psi)}{s^{(N)}_{k+1}-s^{(N)}_k},
	\end{equation}
	noting that, as in \eqref{eq:unrestrictedScaleCasoratian}, the RHS is independent of $k$. For $f:I_N^{(a,b)}\to\mathbb R$, denote by $\bar f:I_N^{[a,b]}\to\mathbb R$ its
	zero extension to the boundary: $\bar f(i)=f(i)$ for $i\in I_N^{(a,b)}$ and
	$\bar f(a)=\bar f(b)=0$.
	
	\begin{proposition}
		\label{prop:killed-resolvent-and-green-kernel}
		Fix $\psi>0$. Then, for every $f:I_N^{(a,b)}\to\mathbb R$ and
		$i\in I_N^{(a,b)}$,
		\begin{align}\label{eq:discreteHuntResolventFormula}
			R_\psi^{N,(a,b)}f(i)
			&=R_\psi^N\bar f(i)
			-U_{a,b}^{+}(i;\psi)R_\psi^N\bar f(a)
			-U_{a,b}^{-}(i;\psi)R_\psi^N\bar f(b).
		\end{align}
		Equivalently, for every $1\leq i,j\leq N-1$, the killed Green kernel satisfies
		\begin{align}
			G_\psi^{N,(a,b)}(i,j)
			&=G_\psi^N(i,j)
			-U_{a,b}^{+}(i;\psi)G_\psi^N(a,j)
			-U_{a,b}^{-}(i;\psi)G_\psi^N(b,j)
			\label{eq:killedGreenFromFullGreen}\\
			&=\frac{
				U_{a,b}^{-}(i\wedge j; \psi)
				U_{a,b}^{+}(i\vee j; \psi)
			}{w_\psi^{N,(a,b)}}.
			\label{eq:killedGreenDirichletForm}
		\end{align}
	\end{proposition}
	\begin{proof}
		Splitting the unrestricted occupation integral at $\tau_{a,b}^N$ and invoking the strong Markov property gives
		\begin{align}
			R_\psi^N\bar f(i) &= R_\psi^{N,(a,b)}f(i) + \mathbb E_i^N\left[\int_{\tau_{a,b}^N}^{\infty}e^{-\psi t}\bar f(X_t^N)\,\dd t\right] = R_\psi^{N,(a,b)}f(i) + \mathbb E_i^N\left[e^{-\psi\tau_{a,b}^N}R_\psi^N\bar f(X_{\tau_{a,b}^N}^N)\right] \nonumber \\
			&= R_\psi^{N,(a,b)}f(i) + U_{a,b}^{+}(i;\psi)R_{\psi}^{N}\bar{f}(a) + U_{a,b}^{-}(i;\psi)R_{\psi}^{N}\bar{f}(b).
			\label{eq:unrestricted-resolvent-split-at-boundary}
		\end{align}
		Rearranging \eqref{eq:unrestricted-resolvent-split-at-boundary} gives \eqref{eq:discreteHuntResolventFormula}, whilst substituting \eqref{eq:resolventDefn} and \eqref{eq:killedMoranResolventGreen} into \eqref{eq:unrestricted-resolvent-split-at-boundary} and comparing coefficients gives \eqref{eq:killedGreenFromFullGreen}. The remaining equality in \eqref{eq:killedGreenDirichletForm} follows by substituting \eqref{eq:GreenFunctionsMoran2} into \eqref{eq:killedGreenFromFullGreen}.
	\end{proof}
	
    We point out that \eqref{eq:killedGreenDirichletForm} is the killed analogue of \eqref{eq:GreenFunctionsMoran2}. 

    \begin{corollary}\label{cor:green}
        For $a<i, j<b$, define
        \begin{equation}
            G_0^{N,(a,b)}(i,j):= \int_0^\infty \frac{\mathbb P_i^N(X_t^N=j,\ t<\tau_{a,b}^N)}{\pi_j}\,\dd t.
        \end{equation} 
        Then 
        \begin{equation}\label{eq:GreenFnCorr1Eq1}
            G_0^{N,(a,b)}(i,j)=\frac{(s^{(N)}_{i\wedge j}-s^{(N)}_a)(s^{(N)}_b-s^{(N)}_{i\vee j})}{s^{(N)}_b-s^{(N)}_a}.
        \end{equation}
        In particular, for any function f on the interior states,
        \begin{equation}\label{eq:GreenFnCorr1Eq2}
            E_i^{N}\left[\int_0^{\tau_{a,b}^{N}}f(X_t^{N})\,dt\right] = \sum_{j=a+1}^{b-1}\pi_jG_0^{N,(a,b)}(i,j)f(j).
        \end{equation}
    \end{corollary}
    
	\subsection{Dirichlet spectrum of the killed Moran generator}
	\label{sec:KilledMoranSpectrum}
    We now focus on the case when $a = 0$ and $b = N$, and for brevity we shall introduce the superscript ${}^{\circ}$ to denote restriction to the interior of $I_N$, for instance $I_N^{\circ} := \{1, \dots, N-1\}$.
    
	Observe that by detailed balance $Q_{N}^{\circ}$ is self-adjoint on $\ell^{2}(I_N^{\circ}, \boldsymbol{\pi})$, and for $f : I_{N}^{\circ} \mapsto \mathbb{R}$ with zero extension $\bar{f}$ to $I_N$, we have that
	\begin{align}\label{eq:chapterEightDiscreteDirichletForm}
		-\langle f,Q_N^\circ f\rangle_{\pi,\circ} = \sum_{i=0}^{N-1} \pi_i\lambda_i \bigl(\bar f(i+1)-\bar f(i)\bigr)^2,
	\end{align}
	where $\langle \cdot, \cdot \rangle_{\boldsymbol{\pi}, \circ}$ is the usual inner product but on $I_N^{\circ}$. Since the RHS of \eqref{eq:chapterEightDiscreteDirichletForm} is strictly positive for $f\not\equiv0$, all eigenvalues of $Q_N^\circ$ are strictly negative. Note that by detailed balance $Q_N^\circ$ is similar to the real symmetric irreducible tridiagonal matrix $D_N Q_{N}^{\circ}D_N^{-1}$ where $D_N := \text{diag}(\sqrt{\pi_1}, \dots, \sqrt{\pi_K})$, and thus all eigenvalues of $Q_N^{\circ}$ are simple. Denote by $\kappa_{j,N}$ the $j^{\textnormal{th}}$ eigenvalue with corresponding eigenvector $\varphi_{j,N}$ such that for $n \in \{1, \dots, N-1\}$
	\begin{align}\label{eq:chapterEightKilledEigenvalues}
		Q_N^\circ\varphi_{n,N} = -\kappa_{n,N}\varphi_{n,N}, \qquad 0<\kappa_{1,N}<\cdots<\kappa_{N-1,N}.
	\end{align}
    Now consider the more general class of functions $\{ \chi_{z}^{N} : z \in \mathbb{C} \}$ where each $\chi_z^N$ is the unique solution to
	\begin{align}\label{eq:chapterEightLeftRecurrenceSolution}
		Q_N\chi_z^N(i)=z\chi_z^N(i), \qquad 1\le i\le N-1, \qquad \chi_z^N(0)=0,\quad \chi_z^N(1)=1.
	\end{align}
	Note that $\chi_{z}^{N}(0) = 0$ imposes a left Dirichlet boundary condition, whilst $\chi_{z}^{N}(1) = 1$ fixes the normalization. To obtain the set of eigenfunctions of the killed generator from the class $\{ \chi_{z}^{N} : z \in \mathbb{C} \}$, we must identify the values of $z$ for which the corresponding solution $\chi_{z}^{N}$ satisfies the corresponding Dirichlet boundary condition at $N$. To this end, we can exploit the formulation in \eqref{eq:chapterEightLeftRecurrenceSolution} to get that $\chi_{z}^{N}(N)$ is a polynomial of degree $N-1$ in $z$. Note in particular that every nonzero eigenvector of $Q_N^\circ$ has a nonzero value at $1$, justifying the normalization in \eqref{eq:chapterEightLeftRecurrenceSolution}. Consequently, the zeros of $\chi_{z}^{N}(N)$ are precisely the eigenvalues $-\kappa_{n,N}$, and we define the corresponding eigenfunctions by
	\begin{align}\label{eq:chapterEightExplicitDiscreteEigenfunction}
		\varphi_{n,N}(i):=\chi_{-\kappa_{n,N}}^N(i), \qquad n=1,\dots,N-1, \qquad i \in I_N^{\circ}.
	\end{align}
	\begin{remark}
		The fundamental solutions $S_{\rho_N(z)}^{\pm}(\cdot)$ provide an alternative representation of the Dirichlet eigenfunctions in \eqref{eq:chapterEightExplicitDiscreteEigenfunction}. Whenever $S_{\rho_N(-\kappa_{n,N})}^{\pm}(\cdot)$ are linearly independent, the normalization $\chi_{-\kappa_{n,N}}^{N}(0) = 0$ and $\chi_{-\kappa_{n,N}}^{N}(1) = 1$ gives
		\begin{align*}
			\varphi_{n,N}(i) = \frac{S_{\rho_N(-\kappa_{n,N})}^{+}(i) - S_{\rho_N(-\kappa_{n,N})}^{+}(0)S_{\rho_N(-\kappa_{n,N})}^{-}(i)}{S_{\rho_N(-\kappa_{n,N})}^{+}(1)-S_{\rho_N(-\kappa_{n,N})}^{+}(0)S_{\rho_N(-\kappa_{n,N})}^{-}(1)}.
		\end{align*}
	\end{remark}
	Finally, self-adjointness of $Q_N^{\circ}$ allows us to re-express the killed transition probabilities and Green's function in terms of the eigenbasis $\{(\kappa_{n,N}, \varphi_{n,N})\}_{n=1}^{N-1}$ to get
    \begin{align*}
		\frac{P_t^{N,\circ}(i,j)}{\pi_j} &= \sum_{n=1}^{N-1} e^{-\kappa_{n,N}t} \frac{\varphi_{n,N}(i)\varphi_{n,N}(j)}{\|\varphi_{n,N}\|_{\pi,\circ}^2}, &
		G_\psi^{N,\circ}(i,j) &= \sum_{n=1}^{N-1} \frac{\varphi_{n,N}(i)\varphi_{n,N}(j)}{(\psi+\kappa_{n,N})\|\varphi_{n,N}\|_{\pi,\circ}^2}.
	\end{align*}
	
	\section{Boundary-to-boundary excursions}\label{sec:Excursions}
	We now exploit the results obtained in Section \ref{sec:KilledMoranModel} to study excursions of the Moran model from the states $0$ and $N$. To this end, we first introduce the appropriate local time through which we will be able to index excursions of the Moran model, and subsequently explicitly construct the corresponding switching excursion point process. In particular we will have a sequence of consecutive ``failed'' excursions which start from and end at the same boundary point, until a ``switching'' excursion occurs where an excursion from one boundary ends at the opposite boundary. We then investigate the boundary killed semigroup and entrance law of the process, before deriving the switching rates for the excursion process. We point out that unlike the diffusion case, we need not consider upper bounds on the mutation rates, since both 0 and $N$ remain accessible regardless of the value of $u$ and $v$.
	
	\subsection{Local time and the excursion space}\label{sec:finite-local-time}
	For $j\in I_N$, define the discrete Markov local time of $X^N$ at $j$ by
	\begin{equation}
		L_t^{N,j} :=\frac{1}{\pi_j}\int_0^{t}\mathbf{1}_{\{X_s^N=j\}}\,ds. \label{eq:lattice-markov-local-time}
	\end{equation}
	We will be especially interested in the two boundary cases $j=0$ and $j=N$, whose local times are denoted by $L^{N,0}$ and $L^{N,N}$, respectively. We retain these two local times separately because, after a boundary-switching excursion, subsequent excursions are indexed by the local time at the newly visited endpoint. For every function $f$ on $\{0,1/N,\ldots,1\}$, we have
	\begin{align}\label{eq:discrete-occupation-local-time}
		\sum_{j=0}^N f\left(\frac{j}{N}\right)L_t^{N,j}\pi_j &=\sum_{j=0}^N f\left(\frac{j}{N}\right)\int_0^{t}\mathbf1_{\{X_s^N=j\}}\,ds = \int_0^{t}f\left(\frac{X_s^N}{N}\right)\,ds. 
	\end{align}
    
	For $b\in\partial_N := \{0,N\}$, let $b^*:=N-b$ denote the opposite endpoint. We shall thus denote an excursion initiated at $b$ and killed on switching by $b\to b^*$, mirroring the notation in \cite[Section~4]{JenkinsEtAlWF}. We define excursion intervals recursively: suppose first that $X_0^N=b_0$ for
	some $b_0\in\partial_N$, and put $D_0=0$ and $B_0:=b_0=X_0^N$.  For
	$n\in\mathbb N$, define
	\begin{align}\label{eq:boundary-excursion}
		G_n&:=\inf\{t>D_{n-1}:X_t^N\notin\partial_N\},  \qquad D_n&:=\inf\{t>G_n:X_t^N\in\partial_N\}, \qquad B_n&:=X_{D_n}^N.
	\end{align}
	Thus $(G_n,D_n)$ is the $n^{\text{th}}$ excursion interval away from the boundary and $B_n\in\partial_N$ records the endpoint at which that excursion is killed. The active boundary immediately before the $n^{\textnormal{th}}$ excursion is $B_{n-1}$, whilst the holding time preceding it is $G_n-D_{n-1}$.
	
	For an excursion $\xi$ initiated at $b\in\partial_N$, let
	\begin{equation}\label{eq:boundary-excursion-entrance-state}
		e_b := \begin{cases}
				1, & b=0,\\
				N-1, & b=N,
			\end{cases}
	\end{equation}
	be its deterministic entrance state, and define its lifetime by
	\begin{equation}\label{eq:marked-boundary-excursion-lifetime}
		\zeta(\xi):=\inf\{t\geq0:\xi(t)=\xi(0)\},
	\end{equation}
	where we adopt the convention that $\inf\varnothing:=\infty$. We distinguish this lifetime from the first boundary-hitting time of the path,
	\begin{equation}\label{eq:boundary-hitting-time-functional}
		\tau_\partial^N(\xi):=\inf\{t\geq0:\xi(t)\in\partial_N\},
	\end{equation}
	noting in particular that the boundary-hitting time $\tau_\partial^N(\xi)$ is finite for both returning and switching excursions. In contrast, $\zeta(\xi)$ is finite only for returning excursions and is infinite for switching excursions. In direct analogy with \cite[Definition~4.1]{JenkinsEtAlWF}, define the path-valued excursion space
	\begin{align}
		\label{eq:marked-boundary-excursion-space}
			\mathcal U_N^b:=\bigl\{\xi\in D([0,\infty),\{0,1,\ldots,N\}):{}& \xi(0)=e_b,\quad 0<\tau_\partial^N(\xi)<\infty,\\
			&\xi(t)\notin\partial_N\text{ for }0\leq t<\tau_\partial^N(\xi),\\
			&\xi(t)=\xi(\tau_\partial^N(\xi))\text{ for }t\geq \tau_\partial^N(\xi)\bigr\}.
	\end{align}
	A returning excursion satisfies $\xi(\tau_\partial^N(\xi))=b$ and has finite lifetime, whereas a switching excursion satisfies $\xi(\tau_\partial^N(\xi))=b^*$ and, because it then remains at $b^*$, has infinite lifetime relative to $b$.
	
	An excursion is said to be \emph{returning} when $\zeta(\xi)<\infty$ and \emph{switching} when $\zeta(\xi)=\infty$.  After a returning excursion, the next excursion is again indexed by the same local time at $b$. After a switching excursion, the active boundary changes, and subsequent excursions are indexed by the local time at the new active boundary. 
	
	For an excursion interval with active boundary $B_{n-1}=b$, define
	\begin{equation}\label{eq:marked-discrete-excursion}
		\xi_n^{N,b}\in\mathcal U_N^b, \qquad \xi_n^{N,b}(t) := \begin{cases}
																X_{G_n+t}^N, & 0\leq t<D_n-G_n,\\
																B_n, & t\geq D_n-G_n.
																\end{cases}
	\end{equation}
	Successive excursions from a boundary point are separated by exponential holding times having finite rate. Hence only finitely many Moran excursions can occur on every bounded interval of physical time, such that they admit the chronological decomposition
	\begin{equation}\label{eq:moran-path-reconstruction}
		X_t^N=\sum_{n\geq1}\left[ B_{n-1}\mathbf 1_{[D_{n-1},G_n)}(t) + \xi_n^{N,B_{n-1}}(t-G_n) \mathbf 1_{[G_n,D_n)}(t)\right].
	\end{equation}
	This contrasts with the Wright--Fisher diffusion, where for any time interval we have uncountably many small excursions accumulating at the boundary and thus the excursion process must be indexed by local time; compare \eqref{eq:moran-path-reconstruction} to \cite[equation~(4.2)]{JenkinsEtAlWF}.
	
	\subsection{Excursion point processes}\label{sec:finite-excursion-processes}
	Suppose now that the process $X^N$ is started from $e_b$, and define the boundary-absorbed path as
	\begin{equation}\label{eq:boundary-absorbed-excursion-path}
		X_t^{N,\partial} := \begin{cases}
								X_t^N, & 0\leq t<\tau_\partial^N,\\
								X_{\tau_\partial^N}^N, & t\geq \tau_\partial^N.
							\end{cases}
	\end{equation}
	For $b\in\partial_N$, let $\mu_N^b$ be the probability measure induced by the law of the process $X^N$ on $\mathcal U_N^b$, namely
	\begin{equation}\label{eq:marked-excursion-probability-law}
		\mu_N^b(A) := \mathbb P_{e_b}^N\left((X_t^{N,\partial})_{t\geq0}\in A\right), \qquad A\in\mathcal B(\mathcal U_N^b).
	\end{equation}
    The process $X^N$ spends an $\operatorname{Exp}(q_b)$ amount of time at $b$, where
	\begin{equation}\label{eq:boundary-departure-rate}
		q_b := \begin{cases}
				\lambda_0, & b=0,\\
				\mu_N, & b=N.
			\end{cases}
	\end{equation}
    By the definition of $L^{N,b}$, the corresponding increment in local-time is therefore $\operatorname{Exp}(\pi_bq_b)$, and independent of the subsequent excursion. Consequently, the \emph{excursion (characteristic) measure} for excursions initiated at $b$ is
	\begin{equation}\label{eq:marked-characteristic-excursion-measure}
		n_N^b := \pi_bq_b\,\mu_N^b = \begin{cases}
										\pi_0\lambda_0\,\mu_N^0, & b=0,\\
										\pi_N\mu_N\,\mu_N^N, & b=N.
										\end{cases}
	\end{equation}
    Given the above, we can define for $t>0$ the \emph{entrance law} associated with excursions initiated at $b$ and observed before they hit either boundary as
	\begin{equation}
		\nu_{N,t}^b(A)
		:=n_N^b
		\left(\xi(t)\in A, t<\tau^N_\partial(\xi)\right),
		\qquad A\subseteq I_N^\circ.
		\label{eq:discrete-excursion-entrance-law}
	\end{equation}
    We point out that whilst \eqref{eq:marked-characteristic-excursion-measure} is a measure on the excursion space $\mathcal{U}_{N}^{b}$, \eqref{eq:discrete-excursion-entrance-law} is a measure on $I^{\circ}_N$. 
    
	For an excursion interval with active boundary $B_{n-1}=b$, its local-time index is given by $\ell_n^{N,b}:=L_{G_n}^{N,b}$. When $B_0=b$, let $J_N^b:=\inf\{n\geq1:B_n\neq b\}$ be the index corresponding to the first switching excursion. 
    
    We can describe the excursions of $X^N$ from $b$ through the point process $\Xi_N^b$ on $\mathbb{R}_{+} \times \mathcal{U}_{N}^{b}$, where if at the local time $\ell$, we have that the process $X^N$ performs an excursion $\xi$ away from $b$ then $(\ell, \xi) \in \Xi_N^b$. As previously mentioned, we make a distinction between returning excursions (which return to $b$ at $\zeta(\xi)$), and those that switch to the opposite boundary through the following partition of $\mathcal{U}_N^b$
	\begin{align}
		\mathcal U_{N,r}^b &:= \{\xi\in\mathcal U_N^b : \eta_b(\xi)<\infty\}, & \mathcal U_{N,s}^b &:= \{\xi\in\mathcal U_N^b:\eta_b(\xi)=\infty\}. \label{eq:switching-excursion-subspace}
	\end{align}
	\begin{theorem}[It\^o \cite{Ito}]\label{prop:marked-discrete-excursion-ppp}
		Fix $b\in\partial_N$ and start the chain from $b$.  Then $\Xi_N^b$ is a Poisson point process on $\mathbb R^+\times\mathcal U_N^b$ with intensity $d\ell\otimes n_N^b$, thus
        \begin{enumerate}
            \item For and $A\subset \mathcal U_N^b$ for which $n_N^b(A)<\infty$, we have that 
            \begin{equation*}
                \mathbb P (\text{no point of }\Xi_N^b \text{ fall in }(0,\ell)\times A)=e^{-\ell n_N^b(A)} .
            \end{equation*}
            \item For disjoint $A_1, \ldots A_k \subset \mathcal U_N^b$ with $n_N^b(A_i)<\infty, \forall i$, if $A=\cup_{i=1}^kA_i$ and $T:=\inf\{\ell >0: \Xi_N^b \in A\}$, then 
            \begin{equation*}
                \mathbb P (\Xi_N^b \in A_i)=\frac{n_N^b(A_i)}{n_N^b(A)}.
            \end{equation*}
        \end{enumerate}
	\end{theorem}

    \begin{remark}
        In view of the above, it follows that the rate at which the excursion processes experiences returning and switching excursions is given by
        \begin{align*}
            n_N^b(\mathcal U_{N,r}^b) &= \pi_bq_b\,\mathbb P_{e_b}^N(X_{\tau_\partial^N}^N = b), & n_N^b(\mathcal U_{N,s}^b) &= \pi_bq_b\,\mathbb P_{e_b}^N(X_{\tau_\partial^N}^N\neq b).
        \end{align*}
    \end{remark}
	
	The sequence of excursions initiated from $b$ terminates with the first switching excursion, which reaches the opposite boundary $b^*$. Once at $b^*$, the strong Markov property implies that subsequent excursions are independent of all preceding excursions, and that the process has law equal to that of the Moran process started from $b^*$. Additionally, we observe that the local time at $b$ remains constant whilst the process is experiencing excursions initiated from $b^*$, until the first switching excursion, at which point the local time at $b$ increases by an $\operatorname{Exp}(\pi_bq_b)$ amount as given in \eqref{eq:boundary-departure-rate}. Thus sequentially concatenating excursions from $b$ indexed by the local time at $b$ with excursions from $b^*$ indexed by the local time at $b^*$ produces an alternating sequence of independent excursion families: excursions from 0 up to the first $0\to N$ switching excursion, followed by excursions from $N$ up to the first $N\to0$ switching excursion, and so on. This is the discrete counterpart of the switching and concatenation mechanism in \cite[Definition~4.5 and equation~(4.2)]{JenkinsEtAlWF}.

    The general Markov entrance-law factorization (see \cite[Chapter~III, Section 3.f]{Blumenthal1992}, \cite[Theorem~5]{RogersGuidedTour}, and \cite[Chapter~VI, Theorem~48.1 and (48.3)]{RogersWilliamsVI8}) allows us to decompose the excursion measure as follows.
    	
	\begin{proposition}\label{prop:discrete-excursion-entrance-law}
		The family $(\nu_{N,t}^b)_{t>0}$ is an entrance law for the killed
		semigroup $(P_t^{N,\circ})_{t\geq0}$: for $s,t>0$,
		\begin{equation}
			\nu_{N,t}^b P_s^{N,\circ}
			=\nu_{N,t+s}^b.
			\label{eq:discrete-entrance-law-property}
		\end{equation}
		Moreover, if $0<t_1<\cdots<t_k$ and
		$A_1,\ldots,A_k\subseteq I_N^\circ$, then
		\begin{align}
			& n_N^b\left(
			\xi(t_r)\in A_r\text{ for }r=1,\ldots,k,
			\ t_k<\tau_\partial(g)
			\right)\nonumber\\
			&\qquad=
			\sum_{x_1\in A_1}\cdots\sum_{x_k\in A_k}
			\nu_{N,t_1}^b(\{x_1\})
			\prod_{r=2}^k
			P_{t_r-t_{r-1}}^{N,\circ}(x_{r-1},\{x_r\}).
			\label{eq:discrete-excursion-fdd-factorization}
		\end{align}
	\end{proposition}

    The decomposition of the excursion measure below is the discrete analogue of \cite[Theorem~4.10]{JenkinsEtAlWF}: it separates the entrance of an excursion from its subsequent propagation according to the process killed at the boundary. 
	
	\subsection{Switching intensities}\label{sec:finite-switching-intensities}
	The following result characterizes explicitly the rate at which switching excursions occur from a boundary $b$.
	\begin{corollary}\label{cor:switchingIntensity}
		For every $N$, the two switching intensities are equal and satisfy
		\begin{equation}\label{eq:switchingIntensities}
			n_N^0\bigl(\mathcal U_{N,s}^0\bigr) = n_N^N\bigl(\mathcal U_{N,s}^N\bigr) = \left(\sum_{k=0}^{N-1}\frac1{\lambda_k\pi_k}
			\right)^{-1}=\frac{1}{s^N_N},
		\end{equation}
        where $s^N_N$ is a scale function given in \eqref{eq:scaleFn}.
	\end{corollary}
	\begin{proof}
		The result follows using \eqref{eq:scaleFn}, \eqref{eq:hittingProbabilitySwitching}, \eqref{eq:marked-characteristic-excursion-measure} and detailed balance.
	\end{proof}
	
	The corresponding calculation for the Wright--Fisher diffusion \cite[Corollary~4.12]{JenkinsEtAlWF} shows that the switching intensities from 0 to 1 and from 1 to 0 coincide, although the intuition behind the result is not clear from the proof. The finite Moran calculations developed here give a direct explanation: this equality results from
    measuring the two directions with their respective boundary local times, factors $\pi_0\lambda_0$ and $\pi_N\mu_N$, together with detailed balance compensate for any potential asymmetry induced by mutation probabilities $u$ and $v$.
	
	\section{Joint path statistics}\label{sec:JointPathStatistics}
	The switching-excursion construction gives a regenerative description of population-wide replacement in the continuous-time Moran model. Starting from 0, successive returning excursions may be interpreted as failed invasions, until the first switching excursion reaches $N$ and completes a successful invasion. The analogous construction then starts from $N$, with the roles of the two types interchanged. By the strong Markov property, each returning excursion regenerates the process at the same boundary, whereas a switching excursion initiates the corresponding sequence of invasion attempts from the opposite boundary.
	
	This representation retains information that is lost when replacement is treated solely as a first-passage event between the two boundaries. In particular, it allows us to retain the sequence of unsuccessful invasions preceding replacement and to characterize jointly their number, heights and durations, together with the duration of the final successful invasion.
	
	To this end, as in Section \ref{sec:KilledMoranModel}, we fix $\psi > 0$, define $\rho_N := \rho_N(\psi)$ as in \eqref{eq:rhoSoln}, and recall the functions $U^{+}_{a,b}(i; \psi)$ and $U^{-}_{a,b}(i; \psi)$ as defined in \eqref{eq:UabMinus} and \eqref{eq:UabPlus}, respectively, with boundary conditions $U^{+}_{a,b}(a; \psi) = U^{-}_{a,b}(b; \psi) = 1$ and $U^{+}_{a,b}(b; \psi) = U^{-}_{a,b}(a; \psi) = 0$. 
	
	\subsection{Single excursion statistics and regenerative identities}\label{sec:RegenIdentities}
	In what follows it will be particularly useful to consider the quantities
	\begin{align}\label{eq:FNSNfunctions}
		F_{N}(\psi) &:= \mathbb{E}_{1}[e^{-\psi \tau^{N}_0}\mathbf1_{\{ \tau_0^N < \tau_N^N \}}] = U^{+}_{0,N}(1; \psi) = \frac{S^{-}_{\rho_N}(N)S^{+}_{\rho_N}(1) - S^{-}_{\rho_N}(1)}{S^{-}_{\rho_N}(N)S^{+}_{\rho_N}(0) - 1}, \nonumber \\ 
		S_{N}(\psi) &:= \mathbb{E}_{1}[e^{-\psi \tau^N_N} \mathbf1_{\{ \tau_N^N < \tau_0^N \}}] = U^{-}_{0,N}(1; \psi) = \frac{S^{-}_{\rho_N}(1)S^{+}_{\rho_N}(0) - S^{+}_{\rho_N}(1)}{S^{-}_{\rho_N}(N)S^{+}_{\rho_N}(0) - 1},
	\end{align}
	using $S^{-}_{\rho_N}(0) = S^{+}_{\rho_N}(N) = 1$. In particular, we observe that $F_{N}(0) + S_{N}(0) = 1$. Define 
	\begin{align}\label{eq:pNExpression}
		p_N := S_{N}(0) = \mathbb{P}_{1}[\tau_N^N < \tau_0^N] = \frac{s_{1}^{(N)} - s_{0}^{(N)}}{s_{N}^{(N)} - s_{0}^{(N)}} = \frac{s_1^{(N)}}{s_{N}^{(N)}},
	\end{align}
	where $s^{(N)}_N$ is a scale function given in \eqref{eq:scaleFn}. Then $p_N$ is the probability that a single invasion is successful, whilst $1 - p_N$ is the probability of an unsuccessful invasion. 
	
	Let $W_0 \sim \operatorname{Exp}(\lambda_0)$ be the holding time at $0$, then we have that
	\begin{equation}\label{eq:B0function}
		B_{0}(\psi) := \E[e^{-\psi W_0}]=\frac{\lambda_{0}}{\lambda_{0}+\psi}.
	\end{equation}
	For $1 \le m \le N-1$, define 
    \begin{equation}\label{eq:max_moran}
        M:=\max_{0\le t\le \tau_{\partial}^N}X_t^{N}
    \end{equation}
    and 
	\begin{align}\label{eq:FNmDefinition}
		F_{N,m}(\psi) := \mathbb{E}_{1}\left[e^{-\psi \tau_0^N} \mathbf1_{\{M=m, \tau_0^N < \tau_N^N\}}\right].
	\end{align}
	On the event $\{M=m, \tau_0^N < \tau_N^N\}$, the chain first reaches $m$ before 0, and after its first visit to $m$, it reaches 0 before $m+1$. Thus the strong Markov property at $\tau_m^N$ gives
	\begin{equation}\label{eq:FNmStrongMarkov}
		F_{N,m}(\psi) = U^{-}_{0,m}(1;\psi)U^{+}_{0,m+1}(m;\psi).
	\end{equation}
	Similarly, for a failed invasion, recall the definition of $\zeta$ in \eqref{eq:marked-boundary-excursion-lifetime} and observe that $\zeta = \tau_{\partial}^N$, allowing us to define
	\begin{equation}
		F_N(\psi,z) := \mathbb{E}_{1}\left[ e^{-\psi\zeta} z^M \mathbf1_{\{\tau_0^N < \tau_N^N\}} \right] = \sum_{m=1}^{N-1}z^mF_{N,m}(\psi).
	\end{equation}
	For a successful invasion we get
	\begin{equation}
		S_N(\psi,z) := \mathbb{E}_{1}\left[ e^{-\psi\zeta} z^M \mathbf1_{\{\tau_N^N < \tau_0^N\}} \right] = z^N S_N(\psi).
	\end{equation}
	
	For $2 \le k \le N-1$, let
	\begin{align}\label{eq:ANkCNkEquations}
		A_{N,k}(\psi)&:=\mathbb{E}_{1}\left[e^{-\psi \tau_0^N} \mathbf1_{\{\tau_0^N < \tau_k^N\}} \right], & C_{N,k}(\psi)&:=\mathbb{E}_{1}\left[e^{-\psi \tau_0^N} \mathbf1_{\{\tau_k^N < \tau_0^N < \tau_N^N\}}\right],
	\end{align}
	and observe that since every failed excursion either returns to 0 before reaching $k$, or reaches $k$ before returning to 0, $A_{N,k}(\psi) + C_{N,k}(\psi) = F_N(\psi)$.
	
	Let $G_N$ be the number of excursion attempts, including the successful attempt, required to reach $N$ after a visit to $0$, set $K_N:=G_N-1$ to be the number of returning excursions, and let $\xi_{1}^{N,0}, \dots, \xi_{K_N+1}^{N,0}$ denote the successive excursion paths. Then the regenerative structure induced by the strong Markov property allows the joint distribution of a number of statistics of the sequence of failed excursions preceding replacement through the following result.
	\begin{proposition}\label{prop:JointStatistics}
		Let $g, h$ be bounded measurable functionals on $\mathcal U_N^0$. Define
		\begin{align}
			\mathcal{F}_N\left[g\right] &:= \int_{\mathcal U_{N,\mathrm{return}}^0} g(\xi)\,\mu_N^0(d\xi), & \mathcal{S}_N\left[h\right] &:= \int_{\mathcal U_{N,\mathrm{switch}}^0} h(\xi)\,\mu_N^0(d\xi),
		\end{align}
		  such that $|\mathcal{F}_N[g]|<1$. Then $\mathcal{F}_N[g]$ and $\mathcal{S}_N[h]$ are the expectations of the corresponding functionals over returning and switching excursions respectively, and we have that
		\begin{equation}\label{eq:JointStatisticsEquation}
			\mathbb{E}_{0}\left[ h(\xi_{K_N+1}^{N,0}) \prod_{r=1}^{K_N}g(\xi_r^{N,0})\right] = \frac{\mathcal{S}_N[h]}{1-\mathcal{F}_N[g]}.
		\end{equation}
	\end{proposition}
	\begin{proof}
		The result follows by the regeneration property of successive excursions, and by observing that conditioning on $K_N$ gives
		\begin{align*}
			\mathbb{E}_{0}\left[ h(\xi_{K_N+1}^{N,0}) \prod_{r=1}^{K_N}g(\xi_r^{N,0})\right] = \mathbb{E}_{0}\left[ \sum_{k=0}^{\infty} h(\xi_{k+1}^{N,0}) \prod_{r=1}^{k}g(\xi_r^{N,0})\right] = \sum_{k=0}^{\infty}\mathcal{F}_{N}\left[g\right]^{k}\mathcal{S}_{N}\left[h\right] = \frac{\mathcal{S}_N[h]}{1-\mathcal{F}_N[g]}.
		\end{align*}
	\end{proof}
	
	\subsection{Replacement times and count-duration statistics}
	From \eqref{eq:FNSNfunctions}, it follows immediately that 
	\begin{align*}
		\mathbb{E}_{1}\left[e^{-\psi \tau_0^N} \mid \tau_0^N < \tau_N^N\right] &= \frac{F_{N}(\psi)}{1-p_N}, & \mathbb{E}_{1}\left[e^{-\psi \tau_N^N} \mid \tau_N^N < \tau_0^N\right] &= \frac{S_{N}(\psi)}{p_N},\\
		\mu_f := \mathbb{E}_{1}\left[ \tau_0^N \mid \tau_0^N < \tau_N^N \right] &= -\frac{F_{N}'(0)}{1-p_N}, & \mu_s := \mathbb{E}_{1}\left[ \tau_N^N \mid \tau_N^N < \tau_0^N\right] &= -\frac{S_{N}'(0)}{p_N},\\
		\sigma_{f}^2 := \operatorname{Var}\left[ \tau_0^N \mid \tau_0^N < \tau_N^N \right] &= \frac{F_{N}''(0)}{1-p_N}-(\mu_f)^2, & \sigma_s^2 := \operatorname{Var}\left[ \tau_N^N \mid \tau_N^N < \tau_0^N\right] &= \frac{S_{N}''(0)}{p_N}-(\mu_s)^2.
	\end{align*}
	\subsubsection{The height and joint height-duration of a failed invasion}
	Using \eqref{eq:FNmStrongMarkov} with $\psi=0$, we get that
	\begin{equation}
		F_{N,m}(0)=\frac{s_{1}^{(N)}}{s_{m}^{(N)}}\frac{s_{m+1}^{(N)}-s_m^{(N)}}{s_{m+1}^{(N)}},
	\end{equation}
	which together with \eqref{eq:FNmDefinition} gives the conditional height law of a failed invasion 
	\begin{equation}
		\mathbb{P}_{1}\left[ M = m \mid \tau_0^N < \tau_N^N \right]=  \frac{s_{1}^{(N)}}{s_{m}^{(N)}}\frac{s_{m+1}^{(N)}-s_{m}^{(N)}}{s_{m+1}^{(N)}}\left(1-\frac{s_{1}^{(N)}}{s_{N}^{(N)}}\right)^{-1},
		\qquad 1\le m \le N-1.
	\end{equation}
	This scale-function characterization of the height of a single excursion is classical for birth--death chains; see, for example, \cite[equations~(39), (42)]{Huillet2013}, where the corresponding excursion-height law is derived for a discrete-time birth--death chain and subsequently applied to Moran models.
	
	The same decomposition above gives the joint Laplace transform for the height and duration a failed invasion as
	\begin{equation}
		\mathbb{E}_{1}\left[e^{-\psi\zeta} \mathbf1_{\{M=m, \tau_0^N < \tau_N^N\}} \right] = U^{-}_{0,m}(1;\psi)U^+_{0,m+1}(m;\psi).
	\end{equation}
	\subsubsection{The total time until a successful invasion}
	Recall that $G_N$ was the total number of invasion attempts until a successful population-wide replacement event including the final successful invasion, and $K_N$ denoted the number of unsuccessful invasions. After each failed excursion the process returns to 0, waits an independent $\operatorname{Exp}(\lambda_0)$ time, and then makes another excursion from 1; hence the successive attempts are iid, with success probability $p_N$. Consequently for $r \in \{ 0, 1, \dots \}$ we have
	\begin{equation}
		\mathbb{P}(K_N=r)=(1-p_N)^r p_N.
	\end{equation}
	Let $T_{0\to N}$ be the total time from a visit to $0$ until the first arrival at $N$, including all boundary holding times and all excursion durations. From \eqref{eq:FNSNfunctions} and \eqref{eq:B0function}, we deduce that a failed attempt contributes $B_{0}(\psi)F_{N}(\psi)$ to $T_{0\to N}$, while the final successful attempt contributes $B_{0}(\psi)S_{N}(\psi)$. Thus conditioning on $K_N$ gives
	\begin{align}\label{eq:TotalTimeUntilSuccess}
		\mathbb{E}_{0}\left[e^{-\psi T_{0\to N}}\right] &= B_{0}(\psi)S_N(\psi)\sum_{r=0}^{\infty}(B_{0}(\psi)F_N(\psi))^r = \frac{B_{0}(\psi)S_{N}(\psi)}{1-B_{0}(\psi)F_{N}(\psi)}.
	\end{align}
	Equivalently, note that $T_{0\to N} \overset{d}{=} \sum_{r=1}^{K_N}(W_0 + \zeta_{r}^{f}) + W_0 + \zeta^{s}$, so setting $g(\xi) = h(\xi) = e^{-\psi \zeta(\xi)}$
	in \eqref{eq:JointStatisticsEquation} and accounting for the fact that we have $K_N+1$ holding times at the boundary, gives \eqref{eq:TotalTimeUntilSuccess}. Exploiting this, we can derive the exact mean switching time
	\begin{equation}\label{eq:ExpectedTimeToSwitch}
		\E_0[T_{0\to N}]=\frac{1}{p_N\lambda_0}+\frac{1-p_N}{p_N}\mu_f+\mu_s,
	\end{equation}
	where we observe that the three terms are respectively the total expected boundary holding time, the expected duration accumulated over failed invasions, and the duration of the final successful invasion.
	
	We can perform a similar calculation for the variance, where we note that the mean and variance of the durations of failed and successful invasions are given by
	\begin{align*}
		m_{f} &= \frac{1}{\lambda_0} + \mu_f, & m_{s} &= \frac{1}{\lambda_0} + \mu_s, & v_{f} &= \frac{1}{\lambda_{0}^{2}} + \sigma_{f}^2, & v_{s} &= \frac{1}{\lambda_{0}^{2}} + \sigma_{s}^2,
	\end{align*}
	such that conditioning on $K_N = r$ gives
	\begin{align*}
		\operatorname{Var}(T_{0\to N}) &= \frac{1-p_N}{p_N}v_{f} + \frac{1-p_N}{p_N^2}m_{f}^2 + v_{s}.
	\end{align*}
	
	\subsubsection{Joint count-duration statistics of a replacement event}
	Recall that $\zeta_1^f,\ldots,\zeta_{K_N}^f$ denote the durations of the failed excursions, excluding boundary holding times, and that $\zeta_N^s$ is the duration of the final successful excursion. Then for $T_N^f := \sum_{r=1}^{K_N}\zeta_r^f$, $u\in[0,1]$ and $\psi_1,\psi_2\ge0$, taking
	\begin{align*}
		g(\xi) &= u e^{-\psi_1 H_{\partial_N}(\xi)}, & h(\xi) &= e^{-\psi_2 H_{\partial_N}(\xi)}
	\end{align*}
	in \eqref{eq:JointStatisticsEquation} gives 
	\begin{equation}
		\mathbb{E}_{0}\left[u^{K_N}e^{-\psi_1 T_N^f}e^{-\psi_2\zeta_N^s}\right] = \frac{S_N(\psi_2)}{1-uF_N(\psi_1)},
	\end{equation}
	from which all mixed moments involving the number of failed attempts and the accumulated duration of failed and successful invasions can be obtained.
	
	\subsection{Extremal statistics of failed invasions}
	\subsubsection{The largest failed invasion height}
	Let $M^{f}_{r}$ be the maximum attained by the $r^{\textnormal{th}}$ failed excursion, and set
	\begin{align*}
		M_N^{\max,f} := \begin{cases}
			0,&K_N=0,\\
			\displaystyle\max_{1\le r\le K_N}M_r^f,&K_N\ge1.
		\end{cases}
	\end{align*}
	Setting
	\begin{align*}
		g(\xi) &= e^{-\psi \zeta(\xi)}\mathbf1_{\left\{\max_{t\geq0}\xi(t)<k\right\}}, & h(\xi) &= e^{-\psi\zeta(\xi)},
	\end{align*}
	we can use \eqref{eq:JointStatisticsEquation} to recover 
	\begin{equation}
		\mathbb{E}_{0}\left[e^{-\psi T_{0\to N}} \mathbf1_{\{M_N^{\max,f}<k\}} \right] = \frac{B_{0,N}(\psi)S_N(\psi)}{1-B_{0,N}(\psi)A_{N,k}(\psi)},
	\end{equation}
	where we recall that $A_{N,k}(\psi)$ is as given in \eqref{eq:ANkCNkEquations}. From this, observe that at $\psi=0$ we have 
	\begin{equation}
		\mathbb{P}(M_N^{\max,f}<k) = \sum_{r=0}^{\infty}A_{N,k}(0)^r p_N = \frac{p_N}{1-A_{N,k}(0)}=\frac{s_k^{(N)}}{s_N^{(N)}},
	\end{equation}
	using the scale function identities \eqref{eq:pNExpression} and $A_{N,k}(0) = 1 - s_{1}^{(N)}/s_{k}^{(N)}$, which in particular gives the full law of $M_N^{\max,f}$
	\begin{equation}
		\mathbb{P}(M_N^{\max,f}=m)=\frac{s_{m+1}^{(N)}-s_m^{(N)}}{s_N^{(N)}}, \qquad 0\le m\le N-1.
	\end{equation}
	
	\subsubsection{The number of failed invasions that reach a threshold}
	Fix $k\in\{1,\ldots,N-1\}$ and define
	\begin{equation}
		N_{N,k}^f:=\sum_{r=1}^{K_N}\mathbf1_{\{M_r^f\ge k\}}.
	\end{equation}
	Recall $C_{N,k}(\psi)$ as defined in \eqref{eq:ANkCNkEquations}, and observe that by the strong Markov property
	\begin{equation}
		C_{N,k}(\psi)=U^-_{0,k}(1;\psi)U^+_{0,N}(k;\psi).
	\end{equation}
    Thus it follows that at $\psi=0$,
	\begin{equation}
		C_{N,k}(0)=\mathbb{P}_{1}\left[ \tau_k^{N} < \tau_0^N \right]\mathbb{P}_{k}\left[ \tau_0^N < \tau_N^N\right] = \frac{s_1^{(N)}}{s_k^{(N)}}\frac{s_N^{(N)}-s_k^{(N)}}{s_N^{(N)}},
	\end{equation}
	and moreover that $A_{N,k}(0)+C_{N,k}(0)=1-p_N$, since every failed excursion either stays below $k$ or reaches $k$. Consequently, conditioning on $ K_N$ gives the probability generating function
	\begin{align*}
		\mathbb{E}_{0}\left[z^{N_{N,k}^f}\right] &= \sum_{r=0}^{\infty}\left(A_{N,k}(0) + zC_{N,k}(0)\right)^{r}p_N = \frac{p_N}{1-A_{N,k}(0)-zC_{N,k}(0)}= \frac{s_k^{(N)}}{s_N^{(N)}-z\bigl(s_N^{(N)}-s_k^{(N)}\bigr)},
	\end{align*}
	where in the last equality we made use of the scale function identities. Thus we deduce that $N_{N,k}^f \sim \operatorname{Geom}_{0}(s_k^{(N)}/s_N^{(N)})$, where $\operatorname{Geom}_{0}(p)$ denotes the geometric law on $\{0,1,2,\ldots\}$ with pmf $(1-p)^r p$. Note further that $\{N_{N,k}^f=0\}=\{M_N^{\max,f}<k\}$, and thus
	\begin{equation}
		\mathbb{P}\left[N_{N,k}^f=0\right] = \frac{s_k^{(N)}}{s_N^{(N)}}.
	\end{equation}
	
	By including the holding times at 0, we get that for $\psi_1,\psi_2\ge0$, and $T_N^{f, \text{tot}} = \sum_{k=1}^{K_N}\left(\zeta_k^{f} + W_{0,k}\right)$, $\zeta^{s, \text{tot}} = \zeta^{s} + W_{0,K_N+1}$
	\begin{equation}
		\mathbb{E}_{0}\left[z^{N_{N,k}^f}e^{-\psi_1 T_N^{f, \text{tot}}}e^{-\psi_2\zeta^{s, \text{tot}}}\right] = \frac{B_{0}(\psi_2)S_N(\psi_2)}
		{1-B_{0}(\psi_1)\bigl(A_{N,k}(\psi_1)+zC_{N,k}(\psi_1)\bigr)},
	\end{equation}
	from which mixed moments involving the threshold count and the two duration variables can be obtained.
	
	\subsection{Additive occupation functionals}
	The approach outlined in Section \ref{sec:RegenIdentities} can also be applied to path functionals that depend on the entire excursion trajectory. As an illustration, let
	\begin{equation}
		A(\xi) := \int_0^{\tau_{\partial}^{N}(\xi)} \xi(t)\,dt,
	\end{equation}
	denote the cumulative allele-count occupation of an excursion. For $\psi, u \ge 0$, define for a returning excursion 
	\begin{equation}
		F_N^{\mathrm{occ}}(\psi,u) := \mathbb{E}_{1}\left[ \exp\left\{ -\psi \tau_0^N-u\int_0^{\tau_0^N}X_t\,dt \right\} \mathbf1_{\{\tau_0^N<\tau_N^N\}} \right].
	\end{equation}
	By the Feynman--Kac formula, $v(i):=\mathbb{E}_{i}[ \exp\{ -\psi \tau_0^N-u\int_0^{\tau_0^N}X_t\,dt \} \mathbf1_{\{\tau_0^N<\tau_N^N\}}]$ is the unique solution of
	\begin{equation}\label{eq:occupation-excursion-feynman-kac}
		\left(\psi+ui-Q_N\right)v(i)=0, \qquad 1\leq i\leq N-1, \qquad v(0) = 1 \qquad v(N) = 0.
	\end{equation}
	Applying the regenerative identity \eqref{eq:JointStatisticsEquation} with
	\begin{align*}
	g(\xi)= \exp\left\{-\psi\tau_\partial^N(\xi)-uA(\xi)\right\}, \qquad h\equiv1,
	\end{align*}
	and denoting $A_N^f := \sum_{r=1}^{K_N}A(\xi_{r}^{N,0}$ gives
	\begin{align*}
		\mathbb{E}_{0}\left[ e^{-\psi T_N^{f}-uA_N^f} \right] &= \sum_{r=0}^{\infty} \mathbb{E}_{0}\left[ e^{-\psi T_N^{f}-uA_N^{f}} \mathbf1_{\{K_N = r\}}\right] = \sum_{r=0}^{\infty}(F_{N}^{\mathrm{occ}}(\psi,u))^{r}p_N = \frac{p_N}{1-F_N^{\mathrm{occ}}(\psi,u)},
	\end{align*}
	reducing joint statistics of additive occupation functionals accumulated over a random number of failed invasions to a single finite-state Feynman--Kac problem.
	
	\subsection{A two-boundary switching process}
	The preceding calculations describe the sequence of invasion attempts leading from a visit to 0 to the first subsequent visit to $N$. The same construction applies after this switch: upon reaching $N$, the process regenerates there and begins a sequence of excursions from $N$, terminating with the first switching excursion that reaches 0. Iterating this construction yields the alternating boundary-switching process. We record the corresponding two-boundary renewal structure below.
	
	Define
	\begin{align*}
		\Phi_{N}^{ab}(\psi) &= \mathbb{E}_{e_{a}}\left[e^{-\psi \tau_{\partial}^{N} \mathbf1_{\{ X_{\tau_{\partial}^{N}} = b \}}}\right], & a,b \in \{0,N\},
	\end{align*}
	where we recall that $e_0 = 1$ and $e_{N} = N-1$. Then $\Phi_N^{0N}(\psi)=S_N(\psi)$ and $\Phi_N^{00}(\psi)=F_N(\psi)$. $\Phi_{N}^{NN}$ and $\Phi_{N}^{N0}$ can be obtained from the same two-boundary hitting formulas starting from $N-1$; equivalently they follow by reflection after interchanging the roles of the mutation parameters $u$ and $v$.
	
	The boundary holding transform matrix $B_N$ (accounting for the holding times of the process at 0 and $N$) and $\Phi$ matrix are
	\begin{align*}
		B_N(\psi) &= \begin{pmatrix}
			\dfrac{\lambda_0}{\lambda_0+\psi} & 0\\[6pt]
			0 & \dfrac{\mu_N}{\mu_N+\psi}
		\end{pmatrix}, & 
		\Phi_N(\psi) &=
		\begin{pmatrix}
			\Phi_N^{00}(\psi)&\Phi_N^{0N}(\psi)\\
			\Phi_N^{N0}(\psi)&\Phi_N^{NN}(\psi)
		\end{pmatrix}.
	\end{align*}
	Then $M_N(\psi):=B_N(\psi)\Phi_N(\psi)$ is the transform kernel for one boundary-to-boundary attempt, including the holding times at the starting boundary. Concatenating attempts gives the matrix geometric series
	\begin{equation}
		M_N+M_N^2+M_N^3+\cdots=(I-M_N)^{-1}M_N,
	\end{equation}
	where the RHS exists whenever the series converges. 
    
	\section{Convergence to the Wright--Fisher diffusion}\label{sec:ConvergenceToDiffusion}
    Under an appropriate rescaling of population size and time, the Moran model converges to the Wright--Fisher diffusion process \cite{Cannings,Etheridge2012}. As killing is not a continuous operation on the path space, this does not directly imply the same results for the killed Moran model, which is the central object in the excursion theory. In this Section we study the convergence of the killed process, and in the process obtain results for the corresponding killed Wright--Fisher diffusion. Since we shall also be interested in the convergence of the Moran excursions to the corresponding Wright--Fisher excursions, throughout this section we shall assume that $0<\theta_1,\theta_2<1$, such that both endpoints 0 and 1 of the limiting Wright--Fisher diffusion are regular, such that the diffusion can enter the interior from either endpoint and can subsequently reach either one \cite[Section~2]{JenkinsEtAlWF}.

    \subsection{Excursion scaling and switching-rate convergence}\label{sec:local-time-boundary-excursions}    
    We begin by defining the appropriately scaled Moran model and the corresponding interior state space through
    \begin{equation}
        Z_t^N:=\frac{X_{N^2t}^N}{N},
        \qquad
        E_N^\circ:=\left\{\frac1N,\ldots,\frac{N-1}{N}\right\}.
        \label{eq:scaled-moran-and-interior}
    \end{equation}
    Given this scaling denote by
    \begin{equation}
        \widehat P_t^{N,\circ}\left(\frac{i}{N},\frac{j}{N}\right)
        =P_{N^2t}^{N,\circ}(i,j),
        \qquad i,j\in I_N^\circ.
        \label{eq:scaled-unscaled-killed-semigroups}
    \end{equation}
    the killed semigroup. The stationary probability measure and the correspondingly normalized scale function are
    \begin{equation}
        m_N:=\sum_{i=0}^N\pi_i\delta_{i/N},
        \qquad \widehat s_N(i/N) :=\frac{s_i^{(N)}}{N^2} =\frac1{N^2}\sum_{k=0}^{i-1}\frac1{\pi_k\lambda_k},
    \qquad 0\le i\le N,
        \label{eq:moran-scale-speed}
    \end{equation}
    where $s_i^{(N)}$ is defined in \eqref{eq:scaleFn}. 
    
    Let $Y$ denote the Wright--Fisher diffusion with generator
    \begin{equation}\label{eq:wf-generator-boundary-comparison}
        \mathcal Lf(x)
        =\frac{1}{2}x(1-x)f''(x)
         +\frac{1}{2}\bigl[\theta_1(1-x)-\theta_2x\bigr]f'(x),
    \end{equation}
    having domain $D(\mathcal L) = C^2([0,1])$. Its normalized speed measure and scale function are given by
    \begin{equation}\label{eq:wf-speed-measure-scale-function}
        m(dx):=\frac{x^{\theta_1-1}(1-x)^{\theta_2-1}}
                      {B(\theta_1,\theta_2)}\,dx, \quad 
        s(x):=2B(\theta_1,\theta_2)\int_0^x u^{-\theta_1}(1-u)^{-\theta_2}\,du.
    \end{equation}

    The speed density and scale function have the following behaviour near the boundaries:
    \begin{equation}
        m(x)\asymp x^{\theta_1-1}, \qquad s(x)\asymp x^{1-\theta_1} \quad\text{as }x\downarrow0,
    \end{equation}
    and
    \begin{equation} 
        m(x)\asymp(1-x)^{\theta_2-1}, \qquad s(1)-s(x)\asymp(1-x)^{1-\theta_2} \quad\text{as }x\uparrow1,
    \end{equation}
    where we write $f\asymp g$ if there exist constants $0<c\leq C<\infty$ such that $cg\leq f\leq Cg$. The convergence results showing how the scale function and the speed measure in \eqref{eq:moran-scale-speed} converge uniformly to \eqref{eq:wf-speed-measure-scale-function} are given in Propositions~\ref{prop:moran-stationary-convergence} and~\ref{prop: moran-scale-convergence} in Appendix~\ref{app:principal-eigenpair}.
    
    Recall the local times $L^{N,j}$ as defined in \eqref{eq:lattice-markov-local-time}, as well as the excursion notation from Section~\ref{sec:finite-local-time}. Set
    \begin{equation}\label{restored:eq:lattice-markov-local-time}
        \widehat L_t^{N,j/N}
        :=N^{-2}L_{N^2t}^{N,j}
        =\frac{1}{\pi_j}\int_0^t
          \mathbf 1_{\{Z_s^N=j/N\}}\,ds,
        \qquad 0\leq j\leq N.
    \end{equation}
    In particular, $\widehat L^{N,0}$ and $\widehat L^{N,1}$ index excursions from the two rescaled endpoints 0 and 1 respectively. The occupation identity~\eqref{eq:discrete-occupation-local-time} now becomes
    \begin{equation}\label{eq:discrete-occupation-local-time-2}
        \sum_{j=0}^N f(j/N)\widehat L_t^{N,j/N}\pi_j
        =\int_0^t f(Z_s^N)\,ds.
    \end{equation}
    
    For comparison, let $L^y$ denote the Wright--Fisher Markov local time with the It\^o--McKean normalization relative to $m$, then for any non-negative Borel function $f$
    \begin{equation}\label{eq:wf-markov-local-time-definition}
        \int_{[0,1]}f(y)L_t^y\,m(dy)
        =\int_0^t f(Y_s)\,ds,
    \end{equation}
    see \cite[Section~4, Remarks~4.2 and 4.3]{JenkinsEtAlWF}. Note that in the diffusion case we need to use the Markov local time since the semimartingale local time at both 0 and 1 vanishes \cite[Proposition~4.4]{JenkinsEtAlWF}. In view of the uniform convergence of the stationary measure to the speed measure from Proposition~\ref{prop:moran-stationary-convergence} in Appendix~\ref{app:principal-eigenpair}, it follows that the discrete normalization in \eqref{eq:discrete-occupation-local-time-2} is compatible with the one in \eqref{eq:wf-markov-local-time-definition}.
    
    Recall the excursion measure $n_N^b$ from \eqref{eq:marked-characteristic-excursion-measure}, and define the rescaling of an excursion $\xi$ as
    \begin{equation}\label{eq:excursion-space-time-rescaling}
        (\mathcal S_N^b \xi)(t):=\frac{\xi(N^2t)}{N},
        \qquad t\geq0.
    \end{equation}
    We shall henceforth refer to $\widehat\xi := \mathcal S_N^b \xi$ as a \emph{rescaled excursion}, and note that one unit of rescaled time corresponds to $N^2$ units of the original Moran time. Set $\widehat{\mathcal U}_N^b:=\mathcal S_N^b(\mathcal U_N^b)$, $\widehat{\mathcal U}_{N,s}^b:=\mathcal S_N^b(\mathcal U_{N,s}^b)$ and $\widehat{\mathcal U}_{N,r}^b:=\mathcal S_N^b(\mathcal U_{N,r}^b)$, and define 
    \begin{equation}
        \widehat\mu_N^b:=\mu_N^b\circ(\mathcal S_N^b)^{-1}.
    \end{equation}
    A rescaled excursion $\widehat \xi\in\widehat{\mathcal U}_N^b$ starts at $e_b/N$ and is absorbed at 0 or 1 after time $\tau_\partial^N(\xi)/N^2$. In particular, we set
    \begin{equation}\label{restored:eq:rescaled-boundary-hitting-time-functional}
        \widehat\tau_\partial(\widehat{\xi}):=\inf\{t\geq0:\widehat\xi(t)\in\{0,1\}\}.
    \end{equation}
    Under the map $(\ell,\xi)\mapsto(\ell/N^2,\mathcal S_N^b \xi)$, the corresponding excursion Poisson point process in Theorem~\ref{prop:marked-discrete-excursion-ppp} takes values in $\mathbb{R}_{+} \times \widehat{\mathcal{U}}_N^b$ and has intensity $d\ell\otimes \widehat n_N^b$, where the characteristic excursion measure $\widehat n_N^b$ is obtained by replacing the local-time index $\ell$ by $\ell/N^2$, namely
    \begin{equation}\label{eq:rescaled-excursion-measures}
        \widehat n_N^b:=N^2\pi_bq_b\,\widehat\mu_N^b.
    \end{equation}
    Define the corresponding entrance law as
    \begin{equation}
        \widehat\nu_{N,t}^b(A)
        :=\widehat n_N^b\bigl(\widehat\xi(t)\in A,\ t<\widehat\tau_\partial(\widehat\xi)\bigr)=N^2\pi_bq_b\,
          \widehat P_t^{N,\circ}(e_b/N,A),
        \qquad t>0,\quad A\subseteq E_N^\circ.
           \label{restored:eq:explicit-discrete-entrance-law}
    \end{equation}
    Proposition~\ref{prop:discrete-excursion-entrance-law} applied to the rescaled process gives the following factorization of the excursion measure through $\widehat\nu_{N,t}^b$ and $\widehat P^{N,\circ}_t$.
    
    For $b\in\{0,1\}$, let 
    \begin{align*}
    \mathcal U_{\mathrm{WF}}^b := \Big\{ \xi \in D([0,\infty),[0,1]) : &\xi(0) = b, 0 < T_{\partial}(\xi) < \infty, \\
    &\xi(t) \notin \{0,1\} \text{ for } 0\le t < T_{\partial}(\xi) \\
    &\xi(t) = \xi(T_{\partial}(\xi)) \text{ for } t \ge T_{\partial}(\xi)\Big\},
    \end{align*}
    with
    \begin{equation}\label{eq:partialT}
        T_\partial(\xi):=\inf\{t>0:\xi(t)\in\{0,1\}\}\in(0,\infty),
    \end{equation}
    such that $\xi(t)\in(0,1)$ for $0<t<T_\partial(\xi)$ and $\xi(t)=\xi(T_\partial(\xi))$ for $t\geq T_\partial(\xi)$. Let $n^b$ denote the characteristic measure for the Wright--Fisher excursions, indexed by the Markov local time $L^b$ normalized relative to the speed measure $m$ as in \eqref{eq:wf-markov-local-time-definition}; see~[12, Section~4]. An excursion is returning if $\xi(T_\partial(\xi))=b$ and switching if $\xi(T_\partial(\xi))=1-b$.
    
    For an excursion $\xi$ starting from $b \in \{0,1\}$, define the functional $M_b$ which retrieves the height of the excursion relative to $b$ by
    \begin{equation}
        M_b(\xi):=\sup_{t\geq0}|\xi(t)-b|.
    \end{equation}
    In particular, for an original Moran excursion $\xi\in\mathcal U_N^0$ and its rescaled path $\widehat\xi:=S_N^0\xi$, we have $M_0(\widehat\xi)=M(\xi)/N$, with $M$ as in \eqref{eq:max_moran}. The rescaled excursion process switches boundary precisely when the height functional equals one. The mass assigned by $\widehat n_N^b$ to infinite-lifetime excursions is the rate, on the local-time scale at $b$, at which the active boundary switches to $b^* = 1-b$. Its asymptotic value follows from the discrete scale function, as shown below.
    
    \begin{proposition}\label{prop:two-sided-switching-intensity}
    For every $k\in\{1,\ldots,N\}$ and every $x\in(0,1]$
    \begin{align}
    \widehat n_N^0(M_0\geq k/N)
    &=\frac{1}{\widehat s_N(k/N)}, \quad
    \widehat n_N^N( M_1\geq k/N)
    =\frac{1}{\widehat s_N(1)-\widehat s_N(1-k/N)}, \label{eq:hatnM}\\
    n^0(M_0\geq x)
    &=\frac{1}{s(x)},\quad \quad \quad \quad\,\,\,
    n^1(M_1\geq x)
    =\frac{1}{s(1)-s(1-x)}. \label{eq:nM}
    \end{align}
    Moreover, whenever $k_N\in\{1,\ldots,N\}$ and $k_N/N\to x\in(0,1]$,
    \begin{equation} \label{eq:excursions_makimum_convergence}
    \widehat n_N^0(M_0\geq k_N/N)
    \longrightarrow n^0(M_0\geq x), \quad
    \widehat n_N^N(M_1\geq k_N/N)
    \longrightarrow n^1(M_1\geq x).
    \end{equation}
    In particular, the two Moran switching intensities coincide and satisfy
    \begin{equation}
    \widehat n_N^0
    \bigl(\widehat{\mathcal U}_{N,s}^0\bigr)
    =
    \widehat n_N^N
    \bigl(\widehat{\mathcal U}_{N,s}^N\bigr)
    =
    \frac{1}{\widehat s_N(1)}
    \end{equation}
    and their common value converges to
    \begin{equation}
    \frac{1}{\widehat s_N(1)}
    \longrightarrow
    n^0(M_0=1)
    =
    n^1(M_1=1)
    =
    \frac{1}{s(1)}
    =
    \frac{1}{2B(\theta_1,\theta_2)B(1-\theta_1,1-\theta_2)}.
    \end{equation}
    \end{proposition}
    
    \begin{proof}
    By the definition of $\widehat n_N^0$, \eqref{eq:moran-scale-speed}, and \eqref{eq:hittingProbabilitySwitching} we get $\widehat n_N^0( M_0\geq k/N)=N^2\pi_0\lambda_0 \frac{s_1^{(N)}}{s_k^{(N)}} =(\widehat s_N(k/N))^{-1}$, where in last equality we used $s_1^{(N)}=(\pi_0\lambda_0)^{-1}$. Similarly, but additionally using detailed balance, we get $\widehat n_N^N( M_1\geq k/N)=(\widehat s_N(1)-\widehat s_N(1-k/N))^{-1}$. 
    
    For the Wright--Fisher diffusion, the Markov property applied once the process hits $x\in(0,1)$ gives
    \begin{equation}
        n^0( M_0=1) = n^0(M_0\geq x) \mathbb P_x(T_1<T_0) = n^0(M_0\geq x)\frac{s(x)}{s(1)},
    \end{equation}
    where $T_y$ denotes the first hitting time of $y$. By~[12, Corollary~4.1], $n^0(M_0=1)=n^1(M_1=1)=s(1)^{-1}$, and thus \eqref{eq:nM} follows.
    
    If $k_N/N\to x\in(0,1]$, Proposition~\ref{prop: moran-scale-convergence} and continuity of $s$ yield \eqref{eq:excursions_makimum_convergence}. Switching excursions are precisely those of height one. Taking $k=N$ gives equality of the Moran switching intensities, and taking $x=1$ gives their diffusion counterparts. 
    \end{proof}

    \subsection{Doob's $h$-transform}
    
    From Proposition \ref{prop:discrete-excursion-entrance-law} it follows that the killed semigroup $(P_t^{N,\circ})_{t\geq0}$ together with the entrance law \eqref{eq:discrete-excursion-entrance-law} completely characterise the excursion measure. Doob's $h$-transform can be used to obtain the semigroup as seen in \cite{GriffithsJenkinsSpano2018} in the case of the Wright-Fisher bridge. By suitably reweighting paths by the positive principal eigenfunction, we can work with conservative process that have no killing terms, and from which the killed process' law can be recovered through an explicit change of measure. For the Wright--Fisher diffusion, the $h$-transformed process is again a Wright--Fisher diffusion, but with mutation parameters which render both endpoints inaccessible. For the Moran model, the transformed process is generally no longer a Moran model, but it remains a conservative birth--death process on the interior states $I^{\circ}_N$.
    
    Let $(\kappa_{1,N},\varphi_{1,N})$ be the principal eigenpair from Section~\ref{sec:KilledMoranSpectrum}, normalized by $\varphi_{1,N}(1)=1$. The Perron--Frobenius theorem gives $\varphi_{1,N}(i)>0$ on $I_N^\circ$. Define
    \begin{equation}\label{eq:chapterNineDiscreteGroundStateGenerator}
        (\widetilde Q_Ng)(i)
        :=\frac{1}{\varphi_{1,N}(i)}
          Q_N^\circ(\varphi_{1,N}g)(i)
          +\kappa_{1,N}g(i),
    \end{equation}
    
    \begin{proposition}\label{prop:finite-ground-state-transform}
    The operator $\widetilde Q_N$ generates a conservative birth--death process on $I_N^\circ$, with rates
    \begin{align}\label{eq:chapterNineTransformedRates}
        \widetilde\lambda_i^N
        &:=\lambda_i
          \frac{\varphi_{1,N}(i+1)}{\varphi_{1,N}(i)},
          &&1\leq i\leq N-2,\nonumber\\
        \widetilde\mu_i^N
        &:=\mu_i
          \frac{\varphi_{1,N}(i-1)}{\varphi_{1,N}(i)},
          &&2\leq i\leq N-1,
    \end{align}
    
    where $\widetilde\mu_1^N=\widetilde\lambda_{N-1}^N=0$. Its semigroup $\widetilde P_t^N=\exp(t\widetilde Q_N)$ satisfies
    \begin{equation}\label{eq:chapterNineDiscreteSemigroupTransform}
        P_t^{N,\circ}(i,j)
        =e^{-\kappa_{1,N}t}
          \frac{\varphi_{1,N}(i)}{\varphi_{1,N}(j)}
          \widetilde P_t^N(i,j),
        \qquad i,j\in I_N^\circ.
    \end{equation}
    
    and the chain is reversible with invariant probability
    \begin{equation}\label{eq:chapterNineTransformedStationaryMeasure}
        \widetilde\pi_i^N
        :=\frac{\varphi_{1,N}(i)^2\pi_i}
          {\displaystyle\sum_{k=1}^{N-1}\varphi_{1,N}(k)^2\pi_k},
        \qquad i\in I_N^\circ.
    \end{equation}
    
    \end{proposition}
    \begin{proof}
    The rates \eqref{eq:chapterNineTransformedRates} and the fact that $\tilde{Q}_N$ defines a conservative birth--death process follows immediately from \eqref{eq:chapterNineDiscreteGroundStateGenerator}, \eqref{eq:operator L^N}, and the fact that $(\kappa_{1,N}, \varphi_{1,N})$ is an eigenpair. If $D$ is the diagonal matrix with $i^{\text{th}}$ diagonal entry $\varphi_{1,N}(i)$, then $Q_N^\circ=D(\widetilde Q_N-\kappa_{1,N}I)D^{-1}$; exponentiation gives \eqref{eq:chapterNineDiscreteSemigroupTransform}. Finally, detailed balance for $\boldsymbol\pi$ implies $\widetilde\pi_i^N\widetilde\lambda_i^N = \widetilde\pi_{i+1}^N\widetilde\mu_{i+1}^N$.
    \end{proof}
    
    Let $\mathcal L^\circ$ denote the generator for the Wright--Fisher diffusion killed upon hitting $\{0,1\}$, with Dirichlet boundary conditions \cite[Proposition~3.2]{JenkinsEtAlWF}. 
    
    Set
    \begin{equation}\label{eq:wright-fisher-ground-state}
        h_{0,1}(x):=x^{1-\theta_1}(1-x)^{1-\theta_2},
        \qquad
        \kappa_{0,1}:=\frac{2-\theta_1-\theta_2}{2}.
    \end{equation}
    
    and observe that $h_{0,1}\in D(\mathcal{L})$. For every compactly supported twice continuously differentiable function $g\in C_c^2(0,1)$, define
    \begin{equation}\label{eq:wright-fisher-transformed-generator-definition}
        \widetilde{\mathcal L}g
        :=h_{0,1}^{-1}\mathcal L^\circ(h_{0,1}g)
          +\kappa_{0,1}g,
    \end{equation}
    
    \begin{proposition}\label{prop:wright-fisher-ground-state-transform}
    The function $h_{0,1}$ satisfies
    \begin{equation}\label{eq:wright-fisher-ground-state-identity}
        \mathcal L h_{0,1}(x)=-\kappa_{0,1}h_{0,1}(x),
        \qquad 0<x<1.
    \end{equation}
    The transformed generator is
    \begin{equation}\label{eq:wright-fisher-transformed-generator}
        \widetilde{\mathcal L}g(x)
        =\frac12x(1-x)g''(x)
          +\frac12\bigl[(2-\theta_1)(1-x)
                         -(2-\theta_2)x\bigr]g'(x).
    \end{equation}
    Let $\widetilde P_t$ be the semigroup associated with $\tilde{\mathcal{L}}$, then for every measurable $f$ such that $f/h_{0,1}$ is bounded,
    \begin{equation}\label{eq:wright-fisher-ground-state-semigroup}
        P_t^\circ f(x)
        =e^{-\kappa_{0,1}t}h_{0,1}(x)
          \widetilde P_t\!\left(\frac{f}{h_{0,1}}\right)(x).
    \end{equation}
    \end{proposition}
    \begin{proof}
    Equations \eqref{eq:wright-fisher-ground-state-identity} and \eqref{eq:wright-fisher-transformed-generator} follow by direct calculation. Since $h_{0,1}$ is positive on $(0,1)$ zero at the boundaries, and satisfies \eqref{eq:wright-fisher-ground-state-identity}, it is an eigenfunction and satisfies $P_t^\circ h_{0,1}=e^{-\kappa_{0,1}t}h_{0,1}$. The Doob $h$-transform of $P_t^{\circ}$ is defined by $\tilde{P}_t g = e^{\kappa_{0,1}t}h_{0,1}^{-1}P_t^{\circ}(h_{0,1}g)$, thus it follows that $\tilde{P}_t$ defines a Markov semigroup whose associated generator coincides with that of a Wright--Fisher diffusion having mutation parameters $2-\theta_1, 2-\theta_2$. 
    \end{proof}

    The following corollary is immediate.
    \begin{corollary}
        Let $p^\circ(t;x,y)$ denote the transition density, with respect to Lebesgue measure, of the Wright--Fisher diffusion with mutation parameters $0<\theta_1,\theta_2<1$, killed upon reaching $\{0,1\}$. Let $\widetilde p(t;x,y)$ denote the transition density of the Wright--Fisher diffusion with mutation parameters $(2-\theta_1,2-\theta_2)$. Then, for $t>0$ and $x,y\in(0,1)$,
        \begin{equation}\label{eq:killedTransitionDensity}
        p^\circ(t;x,y) = e^{-\kappa_{0,1}t} \frac{h_{0,1}(x)}{h_{0,1}(y)} \widetilde p(t;x,y).
        \end{equation}
    \end{corollary}
    The identity \eqref{eq:killedTransitionDensity} expresses the killed Wright--Fisher transition density in terms of that of a Wright--Fisher diffusion with inaccessible boundaries. In particular, it allows the classical spectral and probabilistic representations of the latter to be used to evaluate both the killed transition density, as well as expectations of functionals of the process prior to boundary absorption. The identity also admits an importance-sampling interpretation, since the transformed diffusion can be sampled using existing exact simulation methods \cite{JenkinsSpano,Santetal2026}. A similar result in the context of Wright--Fisher bridges and for the case when at least one of the mutation parameters is 0 can be found in \cite[Section 3.2]{GriffithsJenkinsSpano2018}.

    \begin{remark}
        Note that whilst the $h$-transform of a Wright--Fisher diffusion killed at the boundaries is again a Wright--Fisher diffusion by Proposition \ref{prop:wright-fisher-ground-state-transform}, the $h$-transform a Moran model is not in general a Moran model, since the transition rates obtained in \eqref{eq:chapterNineTransformedRates} in Proposition \ref{prop:finite-ground-state-transform} are no longer of the form \eqref{eq:lambdaMu}. Were this not the case, we would have been able to invoke the classical convergence of Moran models to the Wright--Fisher diffusion under appropriate scaling of space and time. To this end, we need to explicitly show how the $h$-transformed Moran model converges to the $h$-transformed Wright--Fisher diffusion, which is the main aim of Sections \ref{sec:principal-eigenvalue-convergence} and \ref{sec:transformed-process-convergence}.
    \end{remark}
    
    \subsection{Convergence of the principal eigenpair}
    \label{sec:principal-eigenvalue-convergence}
    We first establish convergence of the principal eigenvalue and eigenfunction of the killed Moran generator to those of the killed Wright--Fisher diffusion. Although the eigenfunctions of the Moran generator admit representations in terms of hahn polynomials, establishing their convergence directly would require a detailed asymptotic analysis of these representations together with their corresponding eigenvalues. Instead, we exploit the reversibility of the underlying processes, which allows their principal eigenvalues to be characterised through the associated Dirichlet forms. Comparing the discrete and continuous forms then provides a route to spectral convergence without relying on explicit eigenfunction asymptotics. 
    
    Define $\lVert f\rVert_m^2=\int_0^1 f(x)^2m(dx)$, $\langle f,g\rangle_m:=\int_0^1fg\,dm$, and note that the Dirichlet form associated to the generator $\mathcal{L}^{\circ}$ is given by
    \begin{equation}\label{eq:wright-fisher-dirichlet-form}
        \mathcal E(f):=\frac12\int_0^1x(1-x)|f'(x)|^2\,m(dx),
    \end{equation}
    with domain $\mathcal D_0$ equal to the closure of $C_c^\infty(0,1)$ under the norm
    \begin{equation}\label{eq:wright-fisher-form-norm}
        \lVert f\rVert_{\mathcal E}^2
        :=\lVert f\rVert_m^2+\mathcal E(f).
    \end{equation}
    In particular, for functions $f\in D(\mathcal L^{\circ})$, integration by parts gives the
    equivalent scalar-product representation
    \begin{equation}\label{eq:wright-fisher-dirichlet-form-scalar-product}
        \mathcal E(f)=-\langle f,\mathcal L^\circ f\rangle_m.
    \end{equation}
    Recall that $(\kappa_{0,1}, h_{0,1})$ is an eigenpair, and observe further that since $-\mathcal{L}^{\circ}$ is self-adjoint and non-negative on $L^2(m)$, the Rayleigh--Ritz variational principle identifies the principal eigenvalue as
    \begin{align*}
        \kappa_{\mathrm{pr}} = \inf_{f \in \mathcal{D}_0 \setminus \{0\}} \frac{\mathcal E (f)}{\|f\|_{m}^2}.
    \end{align*}
    To show that $\kappa_{\mathrm{pr}}=\kappa_{0,1}$, take $f\in C_c^\infty(0,1)$ and set $g=f/h_{0,1}$. Integration by parts coupled with the fact that $\mathcal{L}^{\circ}h_{0,1} = (1/2m(x))(x(1-x)m(x)h_{0,1}'(x))'$ gives that
    \begin{equation}\label{eq:wright-fisher-ground-state-dirichlet-form}
     \mathcal E(f)-\kappa_{0,1}\lVert f\rVert_m^2
     =\frac12\int_0^1x(1-x)h_{0,1}(x)^2|g'(x)|^2\,m(dx)\geq0.
    \end{equation}
    Since $C_c^{\infty}(0,1)$ is dense in $\mathcal D_0$ under $\|\cdot\|_m$, the inequality extends to every $f\in\mathcal D_0$. Consequently, $\kappa_{\mathrm{pr}}\geq\kappa_{0,1}$. On the other hand, since $h_{0,1}\in D(\mathcal L^{\circ})$ is an eigenfunction, $\mathcal E(h_{0,1}) = \kappa_{0,1}\|h_{0,1}\|_m^2$, and hence $\kappa_{\mathrm{pr}}\leq \kappa_{0,1}$. Therefore $(\kappa_{0,1}, h_{0,1})$ are the principal eigenpair.
    
    For $u:I_N^\circ\to\mathbb R$, let $\bar u$ be its zero extension to
    $I_N$ and set
    \begin{equation}\label{eq:rescaled-discrete-dirichlet-form}
        \mathcal E_N(u):=N^2\sum_{i=0}^{N-1}\pi_i\lambda_i\bigl(\bar u(i+1)-\bar u(i)\bigr)^2=-N^2\langle u,Q_N^\circ u\rangle_{\pi,\circ},
    \end{equation}
    where the second equality follows from \eqref{eq:chapterEightDiscreteDirichletForm}. Define
    \begin{equation}\label{eq:normalized-principal-eigenvector}
        v_N(i):=\frac{\varphi_{1,N}(i)}
                           {\lVert\varphi_{1,N}\rVert_{\pi,\circ}},
        \qquad 0\leq i\leq N.
    \end{equation}
    such that $\sum_{i=1}^{N-1}v_N(i)^2\pi_i=1$ and $v_N(0)=v_N(N)=0$. Let $V_N\in C([0,1])$ be the piecewise-linear interpolation of $v_N$ on the mesh $i/N$, $0\leq i\leq N$, so $V_N(0)=V_N(1)=0$.

    To show the eigenvalue convergence we will use the Rayleigh--Ritz variational characterization. By sampling smooth, compactly supported functions we establish an upper bound, while relative compactness and lower semicontinuity will give the corresponding lower bound. 
    
    \begin{proposition}\label{prop:principal-eigenvalue-convergence}
    \label{prop:principal-eigenfunction-convergence}
    Suppose that $0<\theta_1,\theta_2<1$.  Then
    \begin{equation}\label{eq:principal-eigenvalue-convergence}
        N^2\kappa_{1,N}\longrightarrow
        \kappa_{0,1}=\frac{2-\theta_1-\theta_2}{2}.
    \end{equation}
    For every compact interval $K\subset(0,1)$,
    \begin{equation}\label{eq:principal-eigenfunction-convergence}
        \sup_{x\in K}\left|
           v_N(\lfloor Nx\rfloor)
           -\frac{h_{0,1}(x)}{\lVert h_{0,1}\rVert_m}
        \right|\longrightarrow0.
    \end{equation}
    \end{proposition}    
    \begin{proof}  
    Recall from \eqref{eq:chapterEightDiscreteDirichletForm} that the killed generator is self-adjoint, so by the Rayleigh--Ritz variational principle the principal eigenvalue $\kappa_{1,N}$ of $-N^2Q_N^\circ$ satisfies
    \begin{equation}\label{eq:discrete-principal-rayleigh-quotient}
        N^2\kappa_{1,N}=\inf_{u\neq0}
        \frac{\mathcal E_N(u)}{\sum_{i=1}^{N-1}u(i)^2\pi_i}.
    \end{equation}
    For any nonzero $f\in C_c^\infty(0,1)$, set $u_N(i)=f(i/N)$. Then equation~\eqref{eq:moran-beta-binomial-density-limit} of Proposition~\ref{prop:moran-stationary-convergence}, Taylor's theorem, and equation~\eqref{eq:interior-energy-coefficient-convergence} of Lemma~\ref{lem:moran-conductance-estimates} give
    \begin{equation}\label{eq:compact-test-rayleigh-convergence}
        \frac{\mathcal E_N(u_N)}{\sum_{i=1}^{N-1}u_N(i)^2\pi_i}
        \longrightarrow\frac{\mathcal E(f)}{\lVert f\rVert_m^2}.
    \end{equation}
    Taking the infimum over $f$ and using $\mathcal{E}(h_{0,1}) = \kappa_{0,1}\|h_{0,1}\|_m^2$ yield
    \begin{equation}\label{eq:principal-eigenvalue-limsup}
        \limsup_{N\to\infty}N^2\kappa_{1,N}\leq\kappa_{0,1}.
    \end{equation}
    For the lower bound, use the normalized positive principal eigenvectors $v_N$ and their interpolations $V_N$ defined above. Since $\mathcal E_N(v_N)=N^2\kappa_{1,N}$, $\mathcal{E}_N(v_N) \le \kappa_{0,1}$ by \eqref{eq:principal-eigenvalue-limsup}.
    
    We first establish compactness on $[0,1]$. Fix a compact interval $K\subset(0,1)$ and a slightly larger compact interior interval $J$. Equations~\eqref{eq:moran-beta-binomial-density-limit} and \eqref{eq:interior-energy-coefficient-convergence} give $\pi_i\asymp N^{-1}$ and $c_{N,i}\asymp N$ on $J$. Since $V_N'=N(\bar v_N(i+1)-\bar v_N(i))$ on each mesh interval, $\mathcal E_N(v_N) \le \kappa_{0,1}$ controls $\lVert V_N'\rVert_{L^2(K)}$, while the normalization controls $\lVert V_N\rVert_{L^2(K)}$. Thus $V_N$ is uniformly bounded in $H^1(K)$, hence uniformly bounded on $K$, and Cauchy--Schwarz gives $|V_N(y)-V_N(x)|\leq C_K|y-x|^{1/2}$ for $x,y\in K$.
    
    To control the endpoints, since $\bar v_N(0)=0$, Cauchy--Schwarz and \eqref{eq:left-right-boundary-resistance-bound} in Lemma~\ref{lem:moran-conductance-estimates} give, for $1\leq i\leq N/2$,
    \begin{equation}\label{eq:left-boundary-cauchy-schwarz}
        v_N(i)^2\leq\mathcal E_N(v_N)\sum_{k=0}^{i-1}\frac1{c_{N,k}}
        \leq C\mathcal E_N(v_N)(i/N)^{1-\theta_1}.
    \end{equation}
    The reflected estimate follows from right bounds in \eqref{eq:left-right-boundary-resistance-bound}. Boundedness of the $\mathcal{E}_N(v_N)$ and linear interpolation therefore give $|V_N(x)|^2\leq C(x+N^{-1})^{1-\theta_1}$ and $|V_N(1-x)|^2\leq C(x+N^{-1})^{1-\theta_2}$ for $0\leq x\leq1/4$. Since $\theta_1,\theta_2<1$, the functions are uniformly small near each endpoint for all sufficiently large $N$. Together with interior equicontinuity, this gives equicontinuity on $[0,1]$: for any prescribed error, first choose small endpoint neighborhoods on which the values are small, then use the interior increment bound on the remaining interval, allowing overlap with those neighborhoods. The finitely many remaining interpolations are continuous and cause no difficulty. The same estimates give uniform boundedness on $[0,1]$.
    Consequently, by the Arzel\`a--Ascoli theorem, the sequence $(V_N)$ of piecewise-linear interpolations is relatively compact in $C([0,1])$, equipped with the uniform norm.
    
    Pick a subsequence such that $\mathcal{E}_{N_j}(v_{N_j}) \rightarrow \liminf_{N\rightarrow\infty}\mathcal{E}_N(v_N)$. Arzel\`a--Ascoli gives a further subsequence converging uniformly to $v\in C([0,1])$, where $v(0)=v(1)=0$ since $V_N(0) = V_N(1) = 0$ for all $N$. Local $H^1$ bounds and a further diagonal extraction give weak convergence in $H^1_{\mathrm{loc}}(0,1)$ to the same $v$. All subsequent limits are taken along this subsequence. Uniform convergence and stationary-measure convergence \eqref{eq:moran-stationary-measure-convergence} imply $\sum_{i=0}^N V_N(i/N)^2\pi_i\to\int_0^1v(x)^2m(x)\,dx$. The sums equal one, so $\lVert v\rVert_m=1$.
    
    On each mesh interval $V_N'=N(\bar v_N(i+1)-\bar v_N(i))$. Consequently, \eqref{eq:interior-energy-coefficient-convergence} in Lemma~\ref{lem:moran-conductance-estimates} and weak lower semicontinuity give, on the chosen subsequence,
    \begin{equation}\label{eq:interior-energy-liminf}
        \frac12\int_\varepsilon^{1-\varepsilon}
           x(1-x)v'(x)^2m(x)\,dx
        \leq\liminf_{N\to\infty}\mathcal E_N(v_N).
    \end{equation}
    Letting $\varepsilon\downarrow0$ yields
    \begin{equation}\label{eq:principal-eigenfunction-compactness}
        \mathcal E(v)\leq\liminf_{N\to\infty}\mathcal E_N(v_N).
    \end{equation}
    Since the RHS of \eqref{eq:principal-eigenfunction-compactness} is finite and $v$ is zero at the boundaries, Lemma~\ref{lem:wright-fisher-form-domain} implies $v\in\mathcal D_0$. The variational lower bound now gives
    \begin{equation}\label{eq:principal-eigenvalue-liminf}
        \liminf_{N\to\infty}N^2\kappa_{1,N}
        \geq\mathcal E(v)\geq\kappa_{0,1}\lVert v\rVert_m^2
        =\kappa_{0,1}.
    \end{equation}
    Together with \eqref{eq:principal-eigenvalue-limsup}, this proves eigenvalue convergence.
    
    To prove the eigenfunction convergence first note that by the compactness established above, every subsequence of $(V_N)$ has a uniformly convergent further subsequence. Its limit $v$ is nonnegative, belongs to $\mathcal D_0$, and has $\lVert v\rVert_m=1$. Lower semicontinuity and eigenvalue convergence give
    \begin{equation}
     \kappa_{0,1}\leq\mathcal E(v)
     \leq\liminf_{N\to\infty}\mathcal E_N(v_N)=\kappa_{0,1}.
    \end{equation}
    Thus $v$ minimizes the limiting Rayleigh quotient. Since $\kappa_{0,1}$ is simple, $v$ must be proportional to $h_{0,1}$. Non-negativity and normalization give $v = h_{0,1} / \|h_{0,1}\|_{m}$. Every subsequence has a further subsequence converging to the same limit, so the entire sequence $(V_N)$ converges uniformly to $v$, and uniform continuity of this limit yields \eqref{eq:principal-eigenfunction-convergence}. 
    \end{proof}
    
    We require the following result to obtain generator convergence. 
    \begin{proposition}\label{prop:principal-eigenvector-ratios}
    For every compact interval $K\subset(0,1)$, with $x_i=i/N$, 
    \begin{align}
        \sup_{x_i\in K}\left|
           N\left(\frac{\varphi_{1,N}(i+1)}{\varphi_{1,N}(i)}-1\right)
           -\frac{h_{0,1}'(x_i)}{h_{0,1}(x_i)}
        \right|&\longrightarrow0,
        \label{eq:principal-forward-ratio}\\
        \sup_{x_i\in K}\left|
           N\left(\frac{\varphi_{1,N}(i-1)}{\varphi_{1,N}(i)}-1\right)
           +\frac{h_{0,1}'(x_i)}{h_{0,1}(x_i)}
        \right|&\longrightarrow0.
        \label{eq:principal-backward-ratio}
    \end{align}
    \end{proposition}
    \begin{proof}
    With $v_N$ normalized as in \eqref{eq:normalized-principal-eigenvector}, put
    \begin{equation}\label{eq:principal-weighted-difference}
     F_N(x_i):=N^2\pi_i\lambda_i\bigl(v_N(i+1)-v_N(i)\bigr),
     \qquad 0\leq i\leq N-1,
    \end{equation}
    and interpolate linearly. The eigenvalue equation and detailed balance give
    \begin{equation}\label{eq:principal-weighted-difference-recurrence}
     F_N(x_i)-F_N(x_{i-1})
     =-N^2\kappa_{1,N}\pi_i v_N(i),\qquad 1\leq i\leq N-1.
    \end{equation}
    On a compact interior interval $J$ containing $K$ in its interior, Propositions~\ref{prop:principal-eigenvalue-convergence} and~\ref{prop:moran-stationary-convergence} bound the right-hand side by $C_J/N$. Thus the interpolations are uniformly Lipschitz on $J$. Moreover, Lemma~\ref{lem:moran-conductance-estimates} and the previously established bound $\mathcal E_N(v_N) \le \kappa_{0,1}$ give
    \begin{equation}\label{eq:principal-weighted-difference-bound}
     \frac1N\sum_{x_i\in J}|F_N(x_i)|^2
     \leq C_J\sum_{x_i\in J}
     \frac{|F_N(x_i)|^2}{N^2\pi_i\lambda_i}
     \leq C_J\mathcal E_N(v_N)\leq \tilde{C}_J.
    \end{equation}
    The mean-square and Lipschitz bounds imply uniform boundedness on $J$. Arzel\`a--Ascoli and diagonal extraction therefore make $(F_N)$ locally uniformly precompact.
    
    For a subsequential limit $F$ and interior $a<b$, telescope
    \eqref{eq:principal-weighted-difference} to obtain
    \begin{equation}
     v_N(\lfloor Nb\rfloor)-v_N(\lfloor Na\rfloor)
     =\frac1N\sum_{i=\lfloor Na\rfloor}^{\lfloor Nb\rfloor-1}
     \frac{F_N(x_i)}{N\pi_i\lambda_i}.
    \end{equation}
    Since $N\pi_i\lambda_i\to x_i(1-x_i)m(x_i)/2$ locally uniformly, Proposition~\ref{prop:principal-eigenfunction-convergence} yields
    \begin{equation}\label{eq:principal-weighted-difference-identification}
     \frac{h_{0,1}(b)-h_{0,1}(a)}{\lVert h_{0,1}\rVert_m}
     =\int_a^b\frac{2F(x)}{x(1-x)m(x)}\,dx.
    \end{equation}
    Differentiating gives
    \begin{equation}
     F(x)=\frac{x(1-x)m(x)h_{0,1}'(x)}{2\lVert h_{0,1}\rVert_m}.
    \end{equation}
    Hence the whole sequence converges locally uniformly. Using $N(v_N(i+1)-v_N(i))=F_N(x_i)/(N\pi_i\lambda_i)$ gives
    \begin{equation}\label{eq:principal-discrete-derivative-convergence}
     \sup_{x_i\in K}\left|
     N\bigl(v_N(i+1)-v_N(i)\bigr)
     -\frac{h_{0,1}'(x_i)}{\lVert h_{0,1}\rVert_m}
     \right|\longrightarrow0.
    \end{equation}
    Division by $v_N(i)$, whose limit is strictly positive on $K$, proves \eqref{eq:principal-forward-ratio}. Applying the same derivative estimate at $x_{i-1}$ on a slightly larger compact interval and using $v_N(i-1)-v_N(i)=-(v_N(i)-v_N(i-1))$ proves \eqref{eq:principal-backward-ratio}. The normalization cancels in both ratios.
    \end{proof}
    
    \subsection{Convergence of the transformed processes}
    \label{sec:transformed-process-convergence}
    Let $\widetilde Z^N$ be Doob’s $h$-transform of the boundary-killed rescaled process $Z^N$, with $h_N(i/N)=\varphi_{1,N}(i)$. Its state space is $E_N^\circ$, and its generator is
    \begin{equation}\label{eq:transformed-rescaled-generator}
        \widetilde{\mathcal A}_N f(x_i)
        :=N^2\widetilde\lambda_i^N\bigl(f(x_i+N^{-1})-f(x_i)\bigr)
          +N^2\widetilde\mu_i^N\bigl(f(x_i-N^{-1})-f(x_i)\bigr).
    \end{equation}
    Here the rates are given by \eqref{eq:chapterNineTransformedRates}, with $\widetilde\mu_1^N=\widetilde\lambda_{N-1}^N=0$.
    
    \begin{proposition}\label{prop:transformed-generator-convergence}
    For $f\in C^2([0,1])$ and every compact interval $K\subset(0,1)$,
    \begin{equation}\label{eq:transformed-generator-convergence}
        \sup_{x_i\in K}\left|
           \widetilde{\mathcal A}_N f(x_i)-\widetilde{\mathcal L}f(x_i)
        \right|\longrightarrow0,
    \end{equation}
    where $\widetilde{\mathcal L}$ is the Wright--Fisher generator with mutation parameters $(2-\theta_1,2-\theta_2)$ in \eqref{eq:wright-fisher-transformed-generator}.
    \end{proposition}
    
    \begin{proof}
    The original rates and Proposition~\ref{prop:principal-eigenvector-ratios} give, locally uniformly,
    \begin{align}
        N(\widetilde\lambda_i^N-\widetilde\mu_i^N)
        &\longrightarrow
          \frac{\theta_1(1-x_i)-\theta_2x_i}{2}
          +x_i(1-x_i)\frac{h_{0,1}'(x_i)}{h_{0,1}(x_i)}
          \nonumber\\
        &=\frac{(2-\theta_1)(1-x_i)-(2-\theta_2)x_i}{2},
          \label{eq:transformed-drift-limit}\\
        \widetilde\lambda_i^N+\widetilde\mu_i^N
        &\longrightarrow x_i(1-x_i).
          \label{eq:transformed-variance-limit}
    \end{align}
    Taylor's formula in \eqref{eq:transformed-rescaled-generator} now proves \eqref{eq:transformed-generator-convergence}; uniform continuity of $f''$ controls the remainder on compact intervals.
    \end{proof}

    Proposition \eqref{prop:transformed-generator-convergence} gives local convergence of generators. We now extend this to $(0,1)$, which suffices since we observe that as $0 < \theta_1, \theta_2 < 1$, the limiting transformed diffusion has mutation rates $2-\theta_1, 2-\theta_2 > 1$, and thus $\{0,1\}$ are unattainable. 
    
    \begin{proposition}\label{prop:transformed-process-convergence}
    Suppose $i_N/N\to x\in(0,1)$ and $\widetilde Z_0^N=i_N/N$. For every $T>0$,
    \begin{equation}\label{eq:transformed-process-convergence}
        (\widetilde Z_t^N)_{0\leq t\leq T}
        \ \Longrightarrow\ (\widetilde Z_t)_{0\leq t\leq T}
        \quad\text{in }D([0,T],[0,1])\text{ with the Skorokhod }J_1\text{ topology},
    \end{equation}
    where $\widetilde Z$ is the Wright--Fisher diffusion with parameters $(2-\theta_1,2-\theta_2)$ started at $x$.
    \end{proposition}
    
    \begin{proof}
    Fix $0<\varepsilon<\min(x,1-x)/2$, and let $\tau_\varepsilon^N$ and $\widetilde\tau_\varepsilon^N$ denote the exit times from $(\varepsilon,1-\varepsilon)$ for $Z^N$ and $\widetilde Z^N$. From equation \eqref{eq:chapterNineDiscreteSemigroupTransform}, it follows that
    \begin{align*}
        M_t := e^{N^2\kappa_{1,N}t}\frac{v_N(NZ_t^N)}{v_N(i_N)}\mathbf1_{\{t < \tau_{\partial}^{N}\}}
    \end{align*}
    is the Radon--Nikodym derivative between the corresponding process laws, and thus $M_t$ is a martingale. Applying optional stopping with the stopping time $\tau_{\varepsilon}^{N} \wedge T$ gives 
    \begin{equation}\label{eq:transformed-stopped-change-of-measure}
        \widetilde{\mathbb P}_{i_N}
           (\widetilde\tau_\varepsilon^N\leq T)
        =\frac{1}{v_N(i_N)}\mathbb E_{i_N}\!\left[
           e^{N^2\kappa_{1,N}\tau_\varepsilon^N}
           v_N(NZ^N_{\tau_\varepsilon^N})
           \mathbf1_{\{\tau_\varepsilon^N\leq T\}}
           \right].
    \end{equation}
    For large $N$, the only two possible exit states are those converging to $\varepsilon$ and $1-\varepsilon$. Proposition~\ref{prop:principal-eigenvalue-convergence} therefore yields
    \begin{equation}\label{eq:transformed-boundary-containment}
        \limsup_{N\to\infty}\widetilde{\mathbb P}_{i_N}
           (\widetilde\tau_\varepsilon^N\leq T)
        \leq e^{\kappa_{0,1}T}
           \frac{\max\{h_{0,1}(\varepsilon),h_{0,1}(1-\varepsilon)\}}
                {h_{0,1}(x)}
        \overset{\varepsilon\rightarrow0}{\longrightarrow}0,
    \end{equation}
    giving compact containment in $(0,1)$. We apply \cite[Chapter~4, Corollary~8.16, with condition~(f) of Corollary~8.7]{EthierKurtz} on $(0,1)$, equipped with any complete metric which generates the same topology as Euclidean metric. For $f\in C_c^\infty(0,1)$, its mesh restriction $f_N$ is uniformly bounded, and Proposition~\ref{prop:transformed-generator-convergence} gives
    \begin{equation}\label{eq:transformed-generator-uniform-test-functions}
        \max_{1\leq i\leq N-1}
        \left|\widetilde{\mathcal A}_N f_N(i/N)
              -\widetilde{\mathcal L}f(i/N)\right|\longrightarrow0.
    \end{equation}
    Outside a compact neighborhood of the support of $f$, both generators vanish for large $N$ because jumps are nearest-neighbor. Thus condition~(f) holds on the whole interior mesh. The limiting martingale problem is well posed by interior smoothness and non-degeneracy, localization, and inaccessibility of both endpoints for parameters $(2-\theta_1,2-\theta_2)$. The cited theorem gives convergence in $D([0,\infty),(0,1))$; continuity of the limit and inclusion into $[0,1]$ give the stated convergence on $[0,T]$.
    \end{proof}

    \subsection{Path-valued convergence of the excursion point processes}
    \label{sec:path-excursion-ppp-convergence}
    
    We now prove convergence of the rescaled excursion point processes introduced in Section~\ref{sec:local-time-boundary-excursions}. We first establish the joint convergence of the stopped excursions, their absorption times, and corresponding terminal endpoints.

    To this end, we introduce the space
    \begin{equation}\label{eq:ppp-mark-space}
     \mathcal X=D([0,\infty),[0,1])\times[0,\infty)\times\{0,1\},
    \end{equation}
    equipped with the local Skorokhod $J_1$, Euclidean and discrete topologies, respectively. We denote by $d_\mathcal{X}$ the product metric on this space. For rescaled Moran excursions, define the map $\Psi_N : \widehat{\mathcal{U}}_{N}^b \mapsto \mathcal X$ such that
    \begin{equation}
        \Psi_N(\xi) :=\bigl(\xi,\widehat\tau_\partial(\xi), \xi(\widehat\tau_\partial(\xi))\bigr), \qquad \xi\in\widehat{\mathcal U}_N^b,\quad b\in\{0,N\}.
    \end{equation}
    For Wright--Fisher excursions, define the map $\Psi : \mathcal U_{\text{WF}}^b \mapsto \mathcal X$ such that
    \begin{equation}
        \Psi(\xi):= \bigl(\xi,T_\partial(\xi),\xi(T_\partial(\xi))\bigr),\qquad  \xi\in\mathcal U_{\mathrm{WF}}^b,\quad b\in\{0,1\}.
    \end{equation}
    We shall refer to the triplets $\Psi_N(\xi)$ and $\Psi(\xi)$ as marks.
    
    Note that the space $\mathcal X$ is Polish as a finite product of Polish spaces. At marks with continuous excursions, convergence in $\mathcal X$ entails convergence of excursion durations and terminal endpoints. Moreover, when the limiting excursion is continuous, convergence in the Skorokhod $J_1$ topology implies uniform convergence of the excursion paths on compact time intervals.
    
    \begin{proposition}
    \label{lem:ppp-stopped-path}
    Recall the definition of $\widehat\tau_{\partial}$ and $T_{\partial}$ in \eqref{restored:eq:rescaled-boundary-hitting-time-functional} and \eqref{eq:partialT} respectively. Denote $T_N=\hat \tau_\partial(Z^N)$ and $T=T_\partial(Y)$. If $Z_0^N=i_N/N\to x=Y_0\in(0,1)$,
    \begin{equation}\label{eq:ppp-stopped-path-limit}
     \bigl((Z^N_{t\wedge T_N})_{t\geq0},T_N,Z^N_{T_N}\bigr)
     \xrightarrow{\mathrm d}
     \bigl((Y_{t\wedge T})_{t\geq0},T,Y_T\bigr), \qquad \text{ in }\mathcal X.
    \end{equation}
    \end{proposition}
    
    \begin{proof}
    We wish to apply the results from Proposition~\ref{prop:transformed-process-convergence} but the transformed processes never hit the boundary. Therefore, fix $0<\varepsilon<\min\{x,1-x\}$, set $a_N=\lfloor N\varepsilon\rfloor$, $b_N=N-a_N$, and observe that for sufficiently large $N$ we have $1 \le a_N < b_N \le N-1$. Define $\tau^N_{\varepsilon}:=\inf\{t\ge0:Z_t^N\in\{a_N/N,b_N/N\}$, and again take $N$ sufficiently large such that the initial state lies strictly between $a_N/N$ and $b_N/N$. Define
    \begin{equation}
        J_{N,\varepsilon}:= \begin{cases}
        0,&Z_{\tau^N_{\varepsilon}}^N=a_N/N,\\
        1,&Z_{\tau^N_{\varepsilon}}^N=b_N/N,
            \end{cases} \qquad
        w_{N,\varepsilon}(t):= \begin{cases}
        Z_t^N,&t<\tau^N_{\varepsilon},\\  
        J_{N,\varepsilon},&t\ge\tau^N_{\varepsilon}.
            \end{cases}
    \end{equation}
    and similarly define $\tau_{\varepsilon},J_\varepsilon,w_\varepsilon$ in an analogous way for $Y$ using $\varepsilon$ and $1-\varepsilon$ instead of $a_N/N$ and $b_N/N$. Denote the corresponding quantities constructed from the Doob's $h$-transformed processes using tildes. For fixed $K>0$ Proposition~\ref{prop:finite-ground-state-transform} implies that for every bounded continuous $F:\mathcal X\to\mathbb R$
    \begin{align}
    &\mathbb E_{i_N}^N
    \left[ F(w_{N,\varepsilon},\tau^N_{\varepsilon},J_{N,\varepsilon})
    \mathbf1_{\{\tau^N_{\varepsilon}\le K\}}
    \right]=\widetilde{\mathbb E}_{i_N}^N
    \left[
    F(\widetilde w_{N,\varepsilon},\widetilde\tau^N_{\varepsilon},\widetilde J_{N,\varepsilon})\mathbf1_{\{\widetilde\tau^N_{\varepsilon}\le K\}} W_{N,\varepsilon}
    \right],
    \end{align}
    where
    \begin{align}
    &W_{N,\varepsilon}:=e^{-N^2\kappa_{1,N}\widetilde\tau^N_{\varepsilon}}
    \frac{v_N(i_N)}
    {v_N(a_N)\mathbf1_{\{\widetilde J_{N,\varepsilon}=0\}}
    +v_N(b_N)\mathbf1_{\{\widetilde J_{N,\varepsilon}=1\}}}.
    \end{align}
    On the other hand, Proposition~\ref{prop:wright-fisher-ground-state-transform} implies that
    \begin{align}
    &\mathbb E_x
    \left[
    F(w_\varepsilon,\tau_\varepsilon,J_\varepsilon)
    \mathbf1_{\{\tau_\varepsilon\le K\}}
    \right]=\widetilde{\mathbb E}_x
    \left[
    F(\widetilde w_\varepsilon,\widetilde\tau_\varepsilon,
    \widetilde J_\varepsilon)
    \mathbf1_{\{\widetilde\tau_\varepsilon\le K\}}W_\varepsilon
    \right],
    \end{align}
    where
    \begin{align}
    &W_\varepsilon:=e^{-\kappa_{0,1}\widetilde\tau_\varepsilon}
    \frac{h_{0,1}(x)}
    {h_{0,1}(\varepsilon)\mathbf1_{\{\widetilde J_\varepsilon=0\}}
    +h_{0,1}(1-\varepsilon)\mathbf1_{\{\widetilde J_\varepsilon=1\}}}.
    \end{align}
    From Proposition~\ref{prop:transformed-process-convergence} it follows that on the event $\{\widetilde\tau_\varepsilon>K\}$, $\tilde Y$ does not exit the interval $(\varepsilon, 1-\varepsilon)$ over the time interval $[0,K]$, so $\tilde{w}_{N,\varepsilon}$ remains inside $(a_N/N, b_N/N)$ over $[0,K]$, for sufficiently large $N$. Moreover, since $\tilde{Y}$ is non-degenerate at both $\varepsilon$ and $1-\varepsilon$, it almost surely crosses the respective level reached at $\tilde{\tau}_{\varepsilon}$ arbitrarily soon after $\tilde{\tau}_{\varepsilon}$. Thus it follows that $\tilde{\tau}_{\varepsilon}^{N} \rightarrow \tilde{\tau}_{\varepsilon}$ almost surely. Consequently, for every bounded continuous $F:\mathcal X\to\mathbb R$,
    \begin{equation}
    \widetilde{\mathbb E}_{i_N}^N
    \left[
    F(\widetilde w_{N,\varepsilon},\widetilde\tau^N_{\varepsilon},\widetilde J_{N,\varepsilon})
    \mathbf1_{\{\widetilde\tau^N_{\varepsilon}\le K\}}
    \right]\longrightarrow
    \widetilde{\mathbb E}_x
    \left[
    F(\widetilde w_\varepsilon,\widetilde\tau_\varepsilon,\widetilde J_\varepsilon)
    \mathbf1_{\{\widetilde\tau_\varepsilon\le K\}}
    \right].
    \end{equation}
    Proposition~\ref{prop:principal-eigenvalue-convergence} gives the convergence of the weights $W_{N,\varepsilon}$ to $W_{\varepsilon},$ and the Skorokhod's representation theorem gives the almost sure convergence of the product, which then by dominated convergence theorem implies
    \begin{equation}
    \mathbb E_{i_N}^N
    \left[
    F(w_{N,\varepsilon},\tau^N_{\varepsilon},J_{N,\varepsilon})
    \mathbf1_{\{\tau^N_{\varepsilon}\le K\}}
    \right]\longrightarrow
    {\mathbb E}_x
    \left[
    F( w_\varepsilon,\tau_\varepsilon, J_\varepsilon)
    \mathbf1_{\{\tau_\varepsilon\le K\}}
    \right].
    \end{equation}
    We thus get the required convergence on $[0,K]$, which we now extend to all of $\mathbb{R}_{+}$. 
    
    From equations \eqref{eq:GreenFnCorr1Eq1} and \eqref{eq:GreenFnCorr1Eq2} of Corollary~\ref{cor:green} it follows after rescaling that $\mathbb E_i^N \left[ T_N\right]\leq \min \{\hat s_N(\frac{i}{N}),\hat s_N(1)-\hat s_N(\frac{i}{N})\}$ for $1\leq i\leq N-1$. Similarly for $Y$ we get $\mathbb E_y \left[ T\right]\leq \min \{ s(y),s(1)-s(y)\}$ for $y\in (0,1)$. Since $\tau^N_{\varepsilon}\leq T_N$ and $\tau_\varepsilon<T$, we can use Markov's inequality to get
    \begin{align}
        \mathbb{E}_{i_N}^{N}\left[F(w_{N,\varepsilon}, \tau_{\varepsilon}^N, J_{N,\varepsilon})\mathbf{1}_{\{\tau_{\varepsilon}^N > K\}}\right] \le \left\|F\right\|_{\infty} \mathbb{P}_{i_N}^{N}\left[\tau_{\varepsilon}^N > K\right] \le \frac{\left\|F\right\|_{\infty}\mathbb{E}_{i_N}^N\left[\tau_{\varepsilon}^N\right]}{K},
    \end{align}
    which therefore yields
    \begin{equation}\label{eq:ppp-inner-stopped-limit}
     (w_{N,\varepsilon},\tau^N_{\varepsilon},J_{N,\varepsilon})
     \xrightarrow{\mathrm d}(w_\varepsilon,\tau_\varepsilon,J_\varepsilon).
    \end{equation}
    It remains to replace the inner stopping levels $\varepsilon$ and $1-\varepsilon$ by the absorbing boundaries. We show that once the process reaches a sufficiently small neighborhood of either boundary, both its remaining lifetime and the probability of returning to a fixed interior level are small. We give the argument near 0, with the case near 1 following through analogous arguments. For $\varepsilon<u<1/2$,
    the scale function formula gives
    \begin{equation}
     \mathbb P^N_{a_N}\left[T_{\lceil Nu\rceil}<T_0\right]
     =\frac{\widehat s_N(a_N/N)}{\widehat s_N(\lceil Nu\rceil/N)},\quad
     \mathbb E_{a_N}\left[T_N\right]\leq\widehat s_N(a_N/N).
    \end{equation}
    Since $\widehat{s}_N \rightarrow s$ uniformly and $s(0) = 0$, both quantities tend to 0 in the limit as $N\rightarrow\infty$ first, and $\varepsilon\rightarrow0$ subsequently. By the strong Markov property at $\tau_{\varepsilon}^N$, these estimates control the process after it reaches either inner level, so $\forall \eta > 0$ we have that $\lim_{\varepsilon\downarrow0}\lim_{N\rightarrow\infty} \mathbb{P}_{i_N}^{N}\left[T_N - \tau_{\varepsilon}^N > \eta\right] = 0$. Thus for any $\delta > 0$ 
    \begin{equation}
    \lim_{\varepsilon\downarrow0}\limsup_{N\to\infty}
    \mathbb P_{i_N}^N
    \bigl(d_{\mathcal X}\bigl((w_{N,\varepsilon},\tau^N_{\varepsilon},J_{N,\varepsilon}),(Z^N_{\cdot\wedge T_N},T_N,Z^N_{T_N})\bigr)>\delta\bigr)=0.    
    \end{equation}
    Combining this with \eqref{eq:ppp-inner-stopped-limit} and using the same argument applied to $Y$ proves \eqref{eq:ppp-stopped-path-limit}.
    \end{proof}

    For $b\in\{0,1\}$, recall $\widehat{n}_N^b$ as defined in \eqref{eq:rescaled-excursion-measures}, and the excursion measure $n^b$ for the Wright--Fisher excursions, and let
    \begin{equation}
    \mathsf n_N^b
    :=\widehat n_N^{b}\circ\Psi_N^{-1},
    \qquad
    \mathsf n^b
    :=n^b\circ\Psi^{-1}
    \end{equation}
    be the images of the rescaled Moran and Wright--Fisher excursion measures under the maps $\Psi_N$ and $\Psi$ into $\mathcal X$. Additionally, let $\Xi_N^b$ and $\Xi^b$ be the auxiliary Poisson point processes on $\mathbb{R}_{+} \times \mathcal X$ with intensities $d\ell\otimes\mathsf n_N^b$ and $d\ell\otimes\mathsf n^b$, respectively.
    
    For $(w,T,J)\in\mathcal X$, write
    \begin{equation}
    M_b(w,T,J):=\sup_{0\leq t\leq T}|w(t)-b|.
    \end{equation}
    For $L<\infty$ and $x\in(0,1)$, the restrictions of $\Xi_N^b$ and $\Xi^b$ to $[0,L]\times\{M_b>x\}$ are finite point measures on $[0,L]\times\mathcal X$, equipped with the weak topology.

    For $b\in\{0,1\}$, define the first switching local-times by
    \begin{align}
        \widehat\sigma_N^b &:=\inf\left\{\ell>0: \Xi_N^b\bigl((0,\ell]\times\{J=1-b\}\bigr)>0 \right\},\\
        \sigma^b &:= \inf\left\{\ell>0: \Xi^b\bigl((0,\ell]\times\{J=1- b\}\bigr)>0 \right\}.
    \end{align}
    where we recall that $J$ denotes the terminal-boundary coordinate of an excursion. Note that by Proposition~\ref{prop:two-sided-switching-intensity}, both local-times $\widehat{\sigma}_N^b$ and $\sigma^b$ are almost surely finite.
    
    \begin{theorem}
    \label{thm:path-excursion-ppp-convergence}
    For every $b\in\{0,1\}$ and $x\in(0,1)$ the finite measures 
    \begin{equation}
        \mathsf n_N^b(\,\cdot\,\cap{\{M_b>x\}}) \overset{\textrm{w}}{\longrightarrow} \mathsf n^b(\,\cdot\,\cap{\{M_b>x\}})
    \end{equation}
    on $\mathcal X$. Moreover,
    \begin{equation}\label{ppp_convergence}
        \left( \widehat\sigma_N^b,\, \left.\Xi_N^b\right|_{[0,\widehat\sigma_N^b]\times \{ M_b>x\}} \right) \overset{\textrm{w}}{\longrightarrow} \left(\sigma^b,\,\left.\Xi^b\right|_{[0,\sigma^b]\times\{M_b>x\}}\right)
    \end{equation}
    in $[0,\infty)\times \mathcal N_f([0,\infty)\times\mathcal X)$, where $\mathcal N_f$ denotes the space of finite point measures equipped with the weak topology.
    \end{theorem}
    \begin{proof}
        Denote $T_N=\hat \tau_\partial(Z^N)$ and $T=T_\partial(Y)$. Fix $0<\varepsilon<\min\{x,1/2\}$ let $a_N=\lceil N\varepsilon\rceil$, $k_N=\lceil Nx \rceil$, $\varepsilon_N=\frac{a_N}{N}$ and define the stopping time $\tau^N_{\varepsilon_N}:=\inf\{t\ge0:Z_t^N=\varepsilon_N\}$. Under $\widehat\mu_N^0(\,\cdot\mid M_0\ge k_N/N)$, we have $\tau^N_{\varepsilon_N}<T_N$, and therefore by the strong Markov property, the shifted $\left((Z^N_{(\tau_{\varepsilon_N}^N+t)\wedge T_N})_{t\ge0},T_N-\tau_{\varepsilon_N}^N,Z^N_{T_N}\right)$ has the same law as $Z^N$ started at $\varepsilon_N$, conditioned on reaching $k_N/N$ before zero. Proposition~\ref{lem:ppp-stopped-path} gives convergence in distribution of the corresponding unconditioned marks. The maximum functional is continuous at marks whose excursions are continuous, and the scale formula gives that $\mathbb P_\varepsilon\!\left(\sup_{t\le T}Y_t\ge v\right)=\frac{s(\varepsilon)}{s(v)}$, for $\varepsilon<v<1$, which implies that the maximum has not atom at $x$, while we observe that $\mathbb{P}_{\varepsilon}(T_x < T_0) > 0$. Consequently, for every bounded continuous $F:\mathcal X\to\mathbb R$,
        \begin{equation}\label{eq:normalised_limit_from_epsilon}
            \frac{ \mathbb E_{a_N}^{N}\!\left[ F\!\left((Z^N_{t\wedge T_N})_{t\ge0},T_N,Z^N_{T_N}\right) \mathbf1_{\{\widehat\tau^N_{k_N/N}<\widehat\tau^N_0\}} \right] }{ \mathbb P_{a_N}^{N}\!\left( \widehat\tau^N_{k_N/N}<\widehat\tau^N_0 \right)}\longrightarrow \frac{ \mathbb E_\varepsilon\!\left[F\!\left((Y_{t\wedge T})_{t\ge0},T,Y_T\right)\mathbf1_{\{T_x<T_0\}}\right]}{\mathbb P_\varepsilon(T_x<T_0)}.
        \end{equation}
        The limit in \eqref{eq:normalised_limit_from_epsilon} describes the diffusion started at $\varepsilon$, conditioned to reach $x$ before zero. As $\varepsilon\downarrow0$, this law converges to the law of a full excursion from zero conditioned to reach $x$, that is $s(x)n^0\bigl(\,\cdot\,\cap{\{M_b>x\}}\bigr)\circ\Psi^{-1}$. It remains to control the corresponding approximation for the Moran excursions. 
        
        Lemma~\ref{lem:conditional_exp_bound} and Markov’s inequality give, for every $h>0$ and all sufficiently large $N$,
        \begin{equation}
            \mathbb E_1^N\left[ \tau_{\varepsilon_N}^N\,\middle|\,\widehat\tau^N_{k_N/N}<\widehat\tau^N_0\right]\le C\varepsilon_N,\qquad\mathbb P_1^N\!\left(\tau_{\varepsilon_N}^N>h\,\middle|\,\widehat\tau^N_{k_N/N}<\widehat\tau^N_0\right)\le\frac{C(\varepsilon+N^{-1})}{h}.
        \end{equation}
        Fix $K>0$ and $h>0$. On the event $\{\tau_{\varepsilon_N}^N\le h\}$, we compare the original stopped excursion with the excursion shifted by $\tau_{\varepsilon_N}^N$ uniformly on $[0,K]$ to get
        \begin{equation}
            \sup_{0\le t\le K}\left|Z^N_{t\wedge T_N}-Z^N_{(\tau_\varepsilon^N+t)\wedge T_N}\right|\le\varepsilon_N+ \sup\limits_{\substack{0\le s,t\le K\\|t-s|\le h}}|Z^N_{(\tau_\varepsilon^N+t)\wedge T_N}-Z^N_{(\tau_\varepsilon^N+s)\wedge T_N}|.
        \end{equation}
    Indeed, for $t < \tau_{\varepsilon_N}^N$ the $Z^N_{t\wedge T_N}$ remains in $[0,\varepsilon_N]$, so \begin{align}
        |Z_{t\wedge T_N}^N - Z_{(\tau_{\varepsilon_N}^N + t) \wedge T_N}^N| \le \varepsilon_N + |Z_{\tau_{\varepsilon_N}\wedge T_N}^N - Z_{(\tau_{\varepsilon_N}^N + t) \wedge T_N}^N|.    
    \end{align}
    For $t \ge \tau_{\varepsilon_N}^N$, the two paths differ only by a time shift of size at most $h$. 
    
    For fixed $\varepsilon$, the shifted marks converge by the first part of the proof. Their limiting laws, in turn, converge as $\varepsilon\downarrow0$ to the continuous-path law $s(x)n^0(\,\cdot\,\cap{\{M_b>x\}})\circ\Psi^{-1}$. Thus for every $\delta>0$,
    \begin{equation}
        \lim_{\varepsilon\downarrow0}\limsup_{N\to\infty} \mathbb P_1^N\!\left( d_{\mathcal X}\!\left( \left((Z^N_{t\wedge T_N})_{t\ge0},T_N,Z^N_{T_N}\right), \left((Z^N_{(\tau_\varepsilon^N+t)\wedge T_N})_{t\ge0}, T_N-\tau_\varepsilon^N,Z^N_{T_N}\right) \right)>\delta \,\middle|\, \widehat\tau^N_{i_N/N}<\widehat\tau^N_0 \right)=0.
    \end{equation}
    Then \cite[Theorem 3.2]{Billingsley1999} gives
    \begin{equation}
        \widehat\mu_N^0(\,\cdot\mid M_0\ge i_N/N)\circ\Psi_N^{-1} \overset{\text{w}}{\longrightarrow} s(x)n^0(\,\cdot\,\cap{\{M_b>x\}})\circ\Psi^{-1}.
    \end{equation}
    Finally, Proposition~\ref{prop:two-sided-switching-intensity} and Proposition~\ref{prop: moran-scale-convergence} imply
    \begin{equation}
        \widehat n_N^0(\,\cdot\,;M_0\ge i_N/N)\circ\Psi_N^{-1}\overset{\text{w}}{\longrightarrow} n^0(\,\cdot\,\cap{\{M_b>x\}})\circ\Psi^{-1}.
    \end{equation}
    We next establish convergence of the excursion point process stopped at their first switching indices. We first restrict our attention to a fixed local-time interval $[0,L]$. Since $\Xi_N^b$ is a Poisson point process with intensity $d\ell \otimes \mathsf n_N^b$, the number of excursions over $[0,L]$ whose height exceeds $x$ is Poisson with mean $L\mathsf n_N^b(M_b>x)$. Conditional on the number of such excursions, the local-time indices of said excursions are independent uniform random variables on $[0,L]$, whilst the excursions themselves are also independent with common law
    \begin{equation}
        \frac{\mathsf n_N^b(\,\cdot\,\cap{\{M_b>x\}})}{\mathsf n_N^b(M_b>x)}.
    \end{equation}
    Therefore
    \begin{equation}
        \Xi_N^b\big|_{[0,L]\times\{M_b>x\}} \overset{\text{w}}{\longrightarrow} \Xi^b\big|_{[0,L]\times\{M_b>x\}}
    \end{equation}
    in the space of finite point measures on $[0,L]\times\mathcal X$, equipped with the weak topology.
    Every switching excursion has $M_b=1$, so the cutoff does not effect the switching points. For a finite point measure on $[0,L]\times\mathcal X$, select the first point having $J=1-b$ and retain all points up to and including it and if none exists, return $L$ and the entire point measure. This map is continuous at finite point measures whose points have distinct local-time indices in $(0,L)$: indeed under weak convergence, the finitely many points can be matched with convergent indices and marks, so their ordering and terminal boundary coordinates are eventually preserved. The limiting Poisson point measure satisfies these conditions almost surely. The continuous mapping theorem therefore gives, for each fixed $L$, joint convergence of the first switching local-time indices truncated at $L$, and the corresponding stopped point measure. On the events $\{\widehat{\sigma}_N^b < L\}$ and $\{\sigma^b < L\}$, these truncated quantities coincide with the untruncated quantities present in \eqref{ppp_convergence}. Finally, Proposition~\ref{prop:two-sided-switching-intensity}, the switching intensities are $\widehat s_N(1)^{-1}$ and $s(1)^{-1}$, respectively. Hence
    \begin{equation}
        \mathbb P(\widehat\sigma_N^b\ge L) =e^{-L/\widehat s_N(1)}, \qquad \mathbb P(\sigma^b\ge L) =e^{-L/s(1)}.
    \end{equation}
    Since $\widehat s_N(1)\to s(1)\in(0,\infty)$,
    \begin{equation}
        \lim_{L\to\infty}\limsup_{N\to\infty}\mathbb P(\widehat\sigma_N^b\ge L)=0, \qquad \lim_{L\to\infty}\mathbb P(\sigma^b\ge L)=0,
    \end{equation}
    from which \eqref{ppp_convergence} follows.
    \end{proof}
    
    \begin{corollary}
    \label{cor:ppp-excursion-statistics}
    Consider the rescaled Moran excursion process from boundary 0, stopped at its first switching excursion to the boundary 1. Let $M_N$ denote the maximum height attained by the returning excursions prior to the first switch, with $M_N = 0$ if there are no returning excursions. Let $C_N(x)$ count the number of returning excursions with height at least $x\in(0,1)$, and let $(W_N^s,T_N^s,J_N^s)$ denote the switching excursion. Then, jointly for any finite collection of thresholds $x\in(0,1)$
    \begin{equation}\label{eq:ppp-statistics-limit}
     (\sigma_N^0,M_N,(C_N(x))_x,W_N^s,T_N^s,J_N^s)
     \overset{\text{w}}{\longrightarrow}
     (\sigma^0,M,(C(x))_x,W^s,T^s,J^s),
    \end{equation}
    where the limiting quantities are defined analogously for the Wright--Fisher excursion process from boundary 0. Their marginal distributions satisfy
    \begin{equation}\label{eq:ppp-height-count-laws}
     \mathbb P(M<x)=\frac{s(x)}{s(1)},\qquad
     C(x)\sim\operatorname{Geom}_0\left(\frac{s(x)}{s(1)}\right).
    \end{equation}
    \end{corollary}
    
    \appendix	
    \section{Auxiliary results}\label{appendix}
    \subsection{Proof of Lemma \ref{lem:gammam}}\label{app:Sec2Proofs}
    \begin{proof}
    	Since
    	\begin{align*}
    		\mathcal{H}_{N,m;i+1}(\ell) &= \mathcal H_{N,m;i}(\ell) \frac{(i+1)(N-i-m+\ell)}{(i+1-\ell)(N-i)}, &
    		\mathcal{H}_{N,m;i-1}(\ell) &= \mathcal{H}_{N,m;i}(\ell)\frac{(i-\ell)(N-i+1)}{i(N-i+1-m+\ell)},
    	\end{align*}
    	it follows that
    	\begin{align*}
    		Q_N \mathcal{H}_{N,m;i}(\ell) &= \frac{\delta_N}{2N^2} (i+a_N) \frac{\ell(N-i) - (m-\ell)(i+1)}{i+1-\ell} \mathcal{H}_{N,m;i}(\ell) \\
    		&{}- \frac{\delta_N}{2N^2} (N-i+b_N) \frac{\ell(N-i+1) - i(m-\ell)}{N-i+1-m+\ell}\mathcal{H}_{N,m;i}(\ell).
    	\end{align*}
    	Observe further that
    	\begin{align}\label{eq:RatioOfMathcalH}
    		\frac{\mathcal{H}_{N,m-1;i}(\ell-1)}{\mathcal{H}_{N,m;i}(\ell)} &= \frac{\ell(N-m+1)}{(i-\ell+1)m}, \nonumber &
    		\frac{\mathcal{H}_{N,m-1;i}(\ell)}{\mathcal{H}_{N,m;i}(\ell)} &= \frac{(m-\ell)(N-m+1)}{m(N-i-m+\ell+1)}.
    	\end{align}
    	which allows us to re-write
    	\begin{align*}
    		Q_N \mathcal{H}_{N,m;i}(\ell) &= -\frac{\delta_N}{2N^2} m(m+a_N+b_N-1) \mathcal{H}_{N,m;i}(\ell) + \frac{\delta_N}{2N^2} m (a_N+\ell-1)\mathcal{H}_{N,m-1;i}(\ell-1) \\
    		&\qquad{}+ \frac{\delta_N}{2N^2} m (m + b_N - \ell -1)\mathcal{H}_{N,m-1;i}(\ell).
    	\end{align*}
    \end{proof}

    \subsection{Technical estimates for excursion convergence}
    \label{app:principal-eigenpair}
    
    The rest of this appendix contains the technical lemmas used in Section~\ref{sec:principal-eigenvalue-convergence}. Throughout, it is assumed that $0<\theta_1,\theta_2<1$, unless otherwise specified.

    The Dirichlet form $\mathcal E$ and its domain $\mathcal D_0$ are defined in Section~\ref{sec:principal-eigenvalue-convergence}, and recall the speed measure introduced in \eqref{eq:wf-speed-measure-scale-function}. We first justify the boundary characterization used there.
    
    \begin{lemma}
    \label{lem:wright-fisher-form-domain}
    Every locally absolutely continuous $f\in L^2(m)$ satisfying $\mathcal E(f)<\infty$ extends continuously to $[0,1]$. If $f(0)=f(1)=0$, then $f\in\mathcal D_0$.
    \end{lemma}
    
    \begin{proof}
    Since $[x(1-x)m(x)]^{-1}\in L^1(0,1)$, Cauchy--Schwarz
    shows that $f'\in L^1(0,1)$, and so $f$ can be extended continuously to $[0,1]$. Suppose now that $f(0)=f(1)=0$. For $0<x<1/2$, again by Cauchy--Schwarz we have
    \begin{equation}\label{eq:wright-fisher-zero-trace-estimate}
        |f(x)|^2
        \leq\left(\int_0^x t(1-t)m(t)|f'(t)|^2\,dt\right)
             \left(\int_0^x\frac{dt}{t(1-t)m(t)}\right)
        \leq Cx^{1-\theta_1}\int_0^x t(1-t)m(t)|f'(t)|^2\,dt.
    \end{equation}
    Choose a smooth cutoff $\chi_\varepsilon$ that is zero on $[0,\varepsilon]$, one on $[2\varepsilon,1]$, takes values in $[0,1]$, and satisfies $|\chi_\varepsilon'|\leq C/\varepsilon$, where $0<\varepsilon<1/4$. Then
    \begin{align}\label{eq:wright-fisher-cutoff-energy}
        &\int_\varepsilon^{2\varepsilon}
           x(1-x)m(x)|f(x)|^2|\chi_\varepsilon'(x)|^2\,dx
        \leq\frac{C}{\varepsilon^2}
           \int_\varepsilon^{2\varepsilon}x^{\theta_1}|f(x)|^2\,dx\nonumber\\
        &\leq\frac{C}{\varepsilon^2}
           \left(\int_0^{2\varepsilon}t(1-t)m(t)|f'(t)|^2\,dt\right)
           \int_\varepsilon^{2\varepsilon}x\,dx
        \leq C\int_0^{2\varepsilon}t(1-t)m(t)|f'(t)|^2\,dt,
    \end{align}
    and the RHS above goes to 0 as $\varepsilon\rightarrow0$. Together with $\int_0^{2\varepsilon}x(1-x)m(x)|f'(x)|^2\,dx\to0$ and $\lVert f-\chi_\varepsilon f\rVert_m\to0$ as $\varepsilon\rightarrow0$, this implies $\lVert f-\chi_\varepsilon f\rVert_{\mathcal E}\to0$. We can replicate the same argument at 1 with $\theta_2$ in place of $\theta_1$. The resulting functions have compact support in $(0,1)$ and can be approximated there by smooth functions in the $\|\cdot\|_{\mathcal{E}}$ norm, because $x(1-x)m(x)$ and $m(x)$ are bounded above and below by positive constants on compact subintervals. Hence $f\in\mathcal D_0$.
    \end{proof}
    
    The next two propositions give the uniform convergence results for the speed measure and the scale functions of the Moran model to the analogous quantities for the Wright--Fisher diffusion.
    
    \begin{proposition}\label{prop:moran-stationary-convergence}
    Let $x \in (0,1)$, and take $x_i=i/N$ such that $\lim_N x_i = x$. For fixed $\theta_1,\theta_2>0$, the stationary measures satisfy
    \begin{equation}\label{eq:moran-stationary-measure-convergence}
     \sum_{i=0}^N\pi_i\delta_{x_i}
     \xrightarrow{\mathrm w}m(x)\,dx.
    \end{equation}
    Uniformly for $0\leq i\leq N$ and all sufficiently large $N$,
    \begin{equation}\label{eq:uniform-beta-binomial-bound}
     \pi_i\asymp\frac{1}{N+1}
     \left(\frac{i+1}{N+1}\right)^{\theta_1-1}
     \left(\frac{N-i+1}{N+1}\right)^{\theta_2-1}.
    \end{equation}
    Moreover, for every $0<\varepsilon<1/2$,
    \begin{equation}\label{eq:moran-beta-binomial-density-limit}
     \max_{\varepsilon\leq i/N\leq1-\varepsilon}
     \left|N\pi_i-m(i/N)\right|\longrightarrow0.
    \end{equation}
    In particular, $N\pi_{i_N}\to m(x)$ whenever $i_N/N\to x\in(0,1)$.
    \end{proposition}
    
    \begin{proof}
    By Proposition~\ref{prop:stat_dis},
    \begin{equation}\label{eq:moran-stationary-gamma-factorization}
     \pi_i=\frac{1}{B(a_N,b_N)}
     \frac{\Gamma(N+1)}{\Gamma(N+a_N+b_N)}
     \frac{\Gamma(i+a_N)}{\Gamma(i+1)}
     \frac{\Gamma(N-i+b_N)}{\Gamma(N-i+1)}.
    \end{equation}
    Stirling's formula gives $\Gamma(k+u)/\Gamma(k+1)\asymp(k+1)^{u-1}$ for $k\geq0$, uniformly for $u$ on compact subsets of $(0,\infty)$. Applying this to the three gamma ratios in \eqref{eq:moran-stationary-gamma-factorization} gives \eqref{eq:uniform-beta-binomial-bound} with exponents $a_N-1,b_N-1$ instead of $\theta_1-, \theta_2-1$. Since $a_N-\theta_1$ and $b_N-\theta_2$ are $O(N^{-1})$ from \eqref{eq:a_N_b_N}, replacing these exponents by $\theta_1-1,\theta_2-1$ changes the expression by a factor tending uniformly to one. On any compact interior interval, the same Stirling estimates give \jaro{We should define $O_\varepsilon(N^{-1})$?}
    \begin{equation}\label{eq:moran-stationary-interior-asymptotic}
     N\pi_i=
     \frac{x_i^{a_N-1}(1-x_i)^{b_N-1}}{B(a_N,b_N)}
     \bigl(1+O_\varepsilon(N^{-1})\bigr),
     \qquad \varepsilon\leq x_i\leq1-\varepsilon,
    \end{equation}
    which proves \eqref{eq:moran-beta-binomial-density-limit}. Finally, the uniform bound \eqref{eq:uniform-beta-binomial-bound} implies
    \begin{equation}\label{eq:moran-stationary-endpoint-mass}
     \sum_{i/N<\varepsilon}\pi_i+
     \sum_{i/N>1-\varepsilon}\pi_i
     \leq C\bigl((\varepsilon+N^{-1})^{\theta_1}
                 +(\varepsilon+N^{-1})^{\theta_2}\bigr),
     \qquad 0<\varepsilon<1/2.
    \end{equation}
    Equation \eqref{eq:moran-stationary-measure-convergence} then follows by first letting $N\to\infty$ and then $\varepsilon\downarrow0$.
    \end{proof}

    The following lemma is particularly useful in establishing convergence of the Moran model Dirichlet forms to those of the Wright--Fisher diffusion.
    \begin{lemma}\label{lem:moran-conductance-estimates}
    Let $c_{N,k}:=N^2\pi_k\lambda_k$ for $0\leq k\leq N-1$. For every compact interval $K\subset(0,1)$,
    \begin{equation}\label{eq:interior-energy-coefficient-convergence}
        \sup_{k/N\in K}\left|
          \frac{c_{N,k}}{N}
          -\frac12\frac{k}{N}\left(1-\frac{k}{N}\right)m(k/N)
        \right|\longrightarrow0.
    \end{equation}
    Uniformly for $0\leq k\leq N/2$,
    \begin{equation}\label{eq:left-boundary-weight-rate-estimates}
        \pi_k\asymp N^{-\theta_1}(k+1)^{\theta_1-1},\qquad
        \lambda_k\asymp\frac{k+1}{N},\qquad
        c_{N,k}\asymp N^{1-\theta_1}(k+1)^{\theta_1}.
    \end{equation}
    At the other endpoint, uniformly for $N/2\leq k\leq N-1$,
    \begin{equation}\label{eq:right-boundary-conductance-estimate}
        c_{N,k}\asymp N^{1-\theta_2}(N-k)^{\theta_2}.
    \end{equation}
    Consequently, for $1\leq i\leq N/2$,
    \begin{equation}
        \sum_{k=0}^{i-1}\frac1{c_{N,k}}
        \leq C(i/N)^{1-\theta_1}, \quad \sum_{k=N-i}^{N-1}\frac1{c_{N,k}}
        \leq C(i/N)^{1-\theta_2}.
        \label{eq:left-right-boundary-resistance-bound}
    \end{equation}
    All comparison constants are independent of $N$ and the indices, for all sufficiently large $N$.
    \end{lemma}
    \begin{proof}
    Proposition~\ref{prop:moran-stationary-convergence} gives locally uniform convergence of $N\pi_k$ to $m(k/N)$. Combining it with \eqref{eq:lambdaMu} gives \eqref{eq:interior-energy-coefficient-convergence}. The global bound \eqref{eq:uniform-beta-binomial-bound} in Proposition \ref{prop:moran-stationary-convergence}, together with $\lambda_k\asymp(N-k)(k+1)/N^2$ (which follows from \eqref{eq:lambdaMu}), give \eqref{eq:left-boundary-weight-rate-estimates} and \eqref{eq:right-boundary-conductance-estimate}. Summing the reciprocal bounds and using $\sum_{j=1}^i j^{-\theta_r}\leq C i^{1-\theta_r}$, for $r=1,2$, gives \eqref{eq:left-right-boundary-resistance-bound}.
    \end{proof}
    
    The following proposition proves the uniform convergence of the Moran model scale functions to those of the Wright--Fisher diffusion.
    \begin{proposition}\label{prop: moran-scale-convergence}
        Recall the scale functions $\widehat s_N$ as in \eqref{eq:moran-scale-speed} and $s$ as in \eqref{eq:wf-speed-measure-scale-function}. Then
        \begin{equation}\label{eq:ppp-uniform-scale}
         \max_{0\leq i\leq N}|\widehat s_N(i/N)-s(i/N)|\longrightarrow0.
        \end{equation}
    \end{proposition}
    \begin{proof}
    Set $c_{N,k}=N^2\pi_k\lambda_k$, then $\widehat s_N(i/N)=\sum_{k=0}^{i-1}\frac{1}{c_{N,k}}$. Since $\frac{1}{s'(x)}$ is bounded away from $0$ on the compact interval $K\subset (0,1)$, equation \eqref{eq:interior-energy-coefficient-convergence} implies that 
    \begin{equation}
        \sum_{k=\lceil Na \rceil}^{\lfloor Nb \rfloor-1}\frac{1}{c_{N,k}}\rightarrow s(b)-s(a)
    \end{equation}
    The bounds \eqref{eq:left-right-boundary-resistance-bound} give
    \begin{equation}
    \sum_{k/N<\varepsilon}\frac1{c_{N,k}}
    +\sum_{k/N>1-\varepsilon}\frac1{c_{N,k}}
    \leq C\left[
    (\varepsilon+N^{-1})^{1-\theta_1}
    +(\varepsilon+N^{-1})^{1-\theta_2}\right].
    \end{equation}
    The integral of $s'$ over the corresponding intervals $(0,\varepsilon)$ and $(1-\varepsilon, 1)$ satisfies a similar bound as the RHS of \eqref{eq:left-right-boundary-resistance-bound} but without the $N^{-1}$ factor. Combining all of the above gives
    \begin{equation}
    \limsup_{N\to\infty}
    \max_{0\leq i\leq N}
    |\widehat s_N(i/N)-s(i/N)|
    \leq C\left[\varepsilon^{1-\theta_1}
    +\varepsilon^{1-\theta_2}\right],
    \end{equation}
    such that taking $\varepsilon\downarrow0$ proves the assertion.
    \end{proof}

    \begin{lemma}\label{lem:conditional_exp_bound}
    There exists a constant $C$, independent of $N,k,\ell$, such that,
    for all sufficiently large $N$ and integers $1\leq \ell<k\leq N$
    with $\ell\leq N/2$,
    \begin{equation}
    \mathbb E_{1}^N
    \left[
    \widehat\tau_{\ell/N}^N
    \,\middle|\,
    \widehat\tau_{k/N}^N<\widehat\tau_{0}^N
    \right]
    \leq C\frac{\ell}{N}.
    \end{equation}
    \end{lemma}
    \begin{proof}
        If $\ell=1$, the hitting time is zero. Suppose $\ell\geq 2$, then
        \begin{align}
            &\mathbb E_1^N \left[ \widehat\tau_{\ell/N}^N \,\middle|\, \widehat\tau_{k/N}^N<\widehat\tau_0^N\right]=\mathbb E_1^N\left[ \widehat\tau_{\ell/N}^N\,\middle|\,\widehat\tau_{\ell/N}^N<\widehat\tau_0^N \right]=\\
            &\frac{1}{\hat s_N(\frac{1}{N})/\hat s_N(\frac{\ell}{N})} \sum_{j=1}^{\ell-1} \pi_j\frac{\hat s_N(\frac{1}{N})(\hat s_N(\frac{\ell}{N})-\hat s_N(\frac{j}{N}))}{\hat s_N(\frac{\ell}{N})}\frac{\hat s_N(\frac{j}{N})}{\hat s_N(\frac{\ell}{N})} \leq \sum_{j=1}^{\ell-1}\pi_j \hat s_N(\frac{j}{N}).
        \end{align}
        Where we have used $\mathbb P_j^N (\widehat\tau_{\ell/N}^N<\widehat\tau_0^N)=\hat s_N(\frac{j}{N})/\hat s_N(\frac{\ell}{N})$ and Corollary~\ref{cor:green}. Lemma~\ref{lem:moran-conductance-estimates}, uniformly for $1\le j\le N/2$ and all sufficiently large $N$ implies that $\pi_j \hat s_N(\frac{j}{N}) \leq C/N$ and summing gives the asserted bound.
    \end{proof} 

    \setcounter{section}{0}
    \setcounter{subsection}{0}
    \gdef\thesection{{}}
    \section{Statement on AI use}
        The authors declare that AI has been used for: Copyediting/proofreading/language editing, Proof assistance, and Proof verification. JS made use of OpenAI's ChatGPT-6 for proof checking and verification, as well as in preparation of the manuscript and presentation. IV made use of OpenAI Codex, using GPT-6 for proof assistance and verification, as well as in formatting and language editing. The authors assume full responsibility for any mistakes.
    \bibliography{draft_bib}
    \bibliographystyle{abbrv}
	
\end{document}